\documentclass{amsart}
\usepackage{amsmath,amssymb,amsthm}
\usepackage{graphics}
\usepackage{boxedeps}
\SetRokickiEPSFSpecial
\HideDisplacementBoxes

\newtheorem{thm}{Theorem}[section]
\newtheorem{lem}{Lemma}[section]
\newtheorem{cor}{Corollary}[section]
\newtheorem{prop}{Proposition}[section]

\newtheorem{exmp}{Example}[section]

\numberwithin{equation}{section}
\def\Z{\Bbb Z}

\def\R{\Bbb R}
\def\P{\Bbb P}
\def\C{\Bbb C}
\def\A{\Bbb A}

\def\flabel#1{\ifmmode #1\else$ #1$\fi}

\let\angle\undefined
\newsymbol\angle 115C
\def\degrees{\ifmmode^\circ\else$^\circ$\fi}

\def \pict #1 by #2 (#3) {\centerline{
\vbox to #2 {\hrule width #1 height 0pt depth 0pt
\vfill{\special{picture #3 }}}}}

\def \picture #1 by #2 (#3 scaled #4) #5{
\dimen0=#1 \dimen1=#2
\divide\dimen0 by 1000 \multiply\dimen0 by #4
\divide\dimen1 by 1000 \multiply\dimen1 by #4
\vbox{\pict \dimen0 by \dimen1 (#3 scaled #4)
\centerline { #5}}}

\catcode`\@=11

\title[]{How tangents solve algebraic equations, \break or a 
remarkable geometry of \break discriminant varieties}

\begin{document}

\author{Gabriel Katz}
\address{Bennington College, Bennington, VT 05201 \hfil\break
\& Department of Mathematics, Brandeis University, Waltham, MA 02454 }

\email{gabrielkatz@rcn.com}

\maketitle

\begin{abstract} Let $\mathcal D_{d,k}$ denote the discriminant variety of degree $d$ 
polynomials in one variable with at least one of its roots being of multiplicity $\geq k$. 
We prove that the tangent cones to $\mathcal D_{d,k}$ span $\mathcal D_{d,k-1}$ thus, 
revealing an extreme ruled nature of these varieties. 
The combinatorics of the web of affine tangent spaces to $\mathcal D_{d,k}$ in 
$\mathcal D_{d,k-1}$ is directly linked to the root multiplicities of the relevant 
polynomials. In fact, solving a polynomial equation $P(z) = 0$ turns out to be equivalent to 
finding hyperplanes through a given point  $P(z)\in \mathcal D_{d,1} \approx \A^d$ 
which are tangent to the discriminant hypersurface $\mathcal D_{d,2}$. We also connect 
the geometry of the Vi\`{e}te map $\mathcal V_d: \A^d_{root} \rightarrow \A^d_{coef}$, given 
by the elementary symmetric polynomials, with the tangents to the discriminant varieties 
$\{\mathcal D_{d,k}\}$.

Various $d$-partitions $\{\mu\}$ provide a refinement $\{\mathcal D_\mu^\circ\}$ of the
stratification  of $\A^d_{coef}$ by the $\mathcal D_{d,k}$'s. 
Our main result, Theorem 7.1, describes an  intricate  
relation between the divisibility of polynomials in one variable and the 
families of spaces tangent to various strata $\{\mathcal D_\mu^\circ\}$.
\end{abstract} 

%\tableofcontents
%\newpage

\bigskip

\section{Introduction}

This exposition depicts a beautiful geometry of stratified discriminant varieties 
which are linked to polynomials in a \emph{single} variable. Perhaps, it was  
Hilbert's ground breaking paper [Hi] which started the exploration. More general discriminant 
varieties have been a focus of an active and broad research 
(cf. [GKZ] which gives a comprehensive account). 
They are studied using methods of algebraic geometry ([GKZ],  [A], [AC], [E], [K]), 
singularity theory ([A1]---[A3], [Va1]---[Va3], [SW], [SK], [GS]) and 
representation theory (with a heavy dose of commutative algebra) ([He], [W1], [W2])
\footnote{This list is far from a complete one: it just reflects some sources 
that I found relevant to this article.}. \smallskip

Here is a text which does not presume an in-depth familiarity with algebraic geometry, 
singularity and representation theories. In fact, it is accessible to a graduate student. 
At the same time, the objects of study are classical and their geometry is fascinating. 
While many basic facts about such discriminants belong to folklore and are spread all over 
the mathematical archipelago, I do not know any self-sufficient elementary treatment 
giving a consistent picture of this small and beautiful island. 

The  discriminants of polynomials in one variable 
constitute a very special class among 
more general discriminants, but it is precisely due to their degenerated nature, 
that they exhibit distinct and unique properties, properties which remain 
uncovered by general theories.     
\bigskip 

This paper has its origins in few  observations that I derived 
from computer-generated  
images of tangent lines to discriminant plane curves (cf. Figures 1 and 8). The flavor 
of the observations can be captured in the slogan: "an algebraic problem of solving  
polynomial equations $$P(z) = z^d + a_1z^{d-1} + ... + a_{d-1}z + a_d = 0$$ is equivalent 
to a geometric problem of finding affine hyperplanes, passing 
through the point $P = (a_1, a_2, \; ... , \; a_d) \in \A^d$ and tangent to the 
discriminant hypersurface $\mathcal D \subset \A^d$" (cf. Corollary 6.1). The discriminant 
hypersurface is comprised of polynomials $P(z)$ with multiple roots, that is, of 
polynomials for which the two equations $\{P(z) = 0,\, P'(z) = 0\}$ have a solution 
$(a_1, a_2, \; ... , \; a_d)$.\footnote{The ground number field is presumed to be $\R$ or $\C$. 
Most of the time, our arguments are not case-sensitive, but their interpretation is.}
 
More generally, one can consider polynomials with roots of multiplicity $\geq k$. 
They form a $(d - k + 1)$-dimensional affine variety $\mathcal D_{d,k} \subset \A^d$.
The resulting stratification $$\A^d = \mathcal D_{d,1} \supset \mathcal D_{d,2} 
\supset  \mathcal D_{d,3} \; ...\; \supset  \mathcal D_{d,d}$$ terminates with 
a smooth curve $\mathcal D_{d,d}$. This stratification has a remarkable property: the 
tangent cones to each stratum $\mathcal D_{d,k}$ span the previous stratum 
$\mathcal D_{d,k-1}$ (Theorem 6.1). Furthermore, $\mathcal D_{d,k-1}$ is comprised 
of the affine subspaces tangent to $\mathcal D_{d,k}$, and the number of 
such subspaces which hit a given point $P \in \mathcal D_{d,k-1}$ is entirely determined 
by the multiplicities of the $P(z)$-roots. 

Surprisingly, the geometry of 
each stratum $\mathcal D_{d,k}$ can be derived from the geometry of a 
single rational curve $\mathcal D_{d,d} \subset \A^d$: its $(d - k + 1)$-st osculating 
spaces span $\mathcal D_{d,k}$ (cf. Theorem 6.2 and [ACGH], pp. 136-137). 
This leads to a "geometrization" 
of the Fundamental Theorem of Algebra (Corollary 6.2). Many of these facts are 
known to experts, but I had a hard time to find out which ones belong to 
folklore, and which ones were actually written down. \smallbreak 

We proceed with a few observations about the $(k - 1)$-dimensional varieties 
$\mathcal D_{d,k}^\vee$ which are the projective \emph{duals} of the varieties $\mathcal D_{d,k}$. 
In Corollary 6.4 we prove that, for $k > 2$,\, 
$deg(\mathcal D_{d,k}^\vee) \leq deg(\mathcal D_{d,k-1})$ (we conjecture that this 
estimate is sharp).\smallskip 
 
In Theorem 6.3 we investigate the interplay between the geometry of the Vi\`{e}te map
$\mathcal V_d: \A^d_{root} \rightarrow \A^d_{coef}$ (given by the elementary 
symmetric polynomials) and the tangents to the discriminant varieties $\{\mathcal D_{d,k}\}$.
\smallskip

Section 7 is devoted to more refined stratification $\{\mathcal D_\mu\}_\mu$ of 
the coefficient space. The strata $\{\mathcal D_\mu\}_\mu$ are indexed by 
$d$-partitions $\{\mu\}$. For a partition $\mu = \{\mu_1 + \mu_2 +  ... + \mu_r = d \}$, 
the variety $\mathcal D_\mu$ is the closure in 
$\A^d_{coef}$ of the set $\mathcal D_\mu^\circ$ of polynomials with $r$ distinct roots 
whose multiplicities are prescribed by the $\mu_i$'s.  
When $\mu = \{k + 1 + 1 + ... + 1 = d\}$, 
$\mathcal D_\mu = \mathcal D_{d,k}$. However, a generic variety $\mathcal D_\mu$ 
exhibits geometric properties very different from the ones of its ruled relative 
$\mathcal D_{d,k}$.\smallskip
 
Our main result is Theorem 7.1. It describes an interesting and intricate  
relation between the divisibility of polynomials in one variable and the 
families of spaces tangent to various strata $\mathcal D_\mu^\circ$'s. 
Among other things, Theorem 7.1 depicts the decomposition of the 
quasiaffine variety $T\mathcal D_\mu^\circ$, comprised of spaces tangent to $\mathcal D_\mu^\circ$, 
into various pieces $\{\mathcal D_{\mu'}^\circ\}$. Also, it is preoccupied 
with the multiplicities of the tangent web forming  $T\mathcal D_\mu^\circ$ (see also Corollary 7.1).
Corollary 7.2 describes a remarkable stabilization of tangent spaces $T_Q\mathcal D_\mu^\circ$, 
as a point $Q \in \mathcal D_\mu^\circ$ approaches one of the singularities 
$\mathcal D_\nu^\circ \subset \mathcal D_\mu$. 
\smallskip

We conclude with a few well-known remarks about the topology of the strata 
$\{\mathcal D_{d,k}^\circ := 
\mathcal D_{d,k} \setminus \mathcal D_{d,k+1}\}$ and $\{\mathcal D_\mu^\circ\}$ in connection 
to the colored braid groups.
\bigskip

After describing the observations above in a draft, I decided that it is a good time 
to consult with experts. I am grateful to Boris Shapiro for an eye-opening education. 
Also, with the help of Harry Tamvakis and Jersey Weyman I have learned about 
a flourishing research  which tackles much more general discriminant 
varieties. My thanks extend to all these people.
\bigbreak 

A perceptive reader might wonder why all our references point towards 
Sections 6 and 7, and what is going on in the other sections of the article. 
The paper is written 
to satisfy two types of readership.  The readers who are willing to endure the pain of 
combinatorics and multiple indices can proceed directly to Section 6, devoted to 
polynomials of a general degree $d$. The readers who prefer to see basic examples 
and special cases ($d = 2, 3$) of theorems from  Sections 6 and 7, 
being striped of combinatorial 
complexities, could be satisfied by the slow pace of Sections 2---5.  In any case, our 
methods are quite elementary and the proofs are self-contained.\smallskip

Some of the graphical images were produced using the \emph{3D-FilmStrip}---a Mac-based 
software tool for a dynamic stereo visualization in geometry. It is developed by 
Richard Palais to whom I am thankful for  help and pleasant conversations. 
  
%%%%%%%%%%%%%%%%%%%%%%%%%%%%%%%%%%

\section{quadratic discriminants}

This section describes some "well-known" and some "less-well-known" geometry of the 
quadratic discriminant. 
It will provide us with a "baby model" of more general geometric structures to come. 
\smallskip

Let $u$ and $v$ be the roots of a monic quadratic polynomial $P(z) = z^2 + bz + c$.  The 
Vi\`{e}te formulas $b = -u - v, \; c = uv$ give rise to a quadratic polynomial map 
$$\mathcal V : (u, v) \rightarrow (-u - v, \, uv)$$ 
from the $uv$-\emph{root plane} $\A^2_{root}$ to the $bc$-\emph{coefficient plane} $\A^2_{coef}$. 
We call it the \emph{Vi\`{e}te Map}. Points of 
the root plane are \emph{ordered} pairs of roots. Therefore, generically, $\mathcal V$ 
is a 2-to-1 map: pairs $(u, v)$ and $(v, u)$ generate the same quadratic 
polynomial. Being restricted to  the diagonal line $L = \{u = v\}$, the map $\mathcal V$ is 
1-to-1.\smallskip

\begin{figure}[ht]
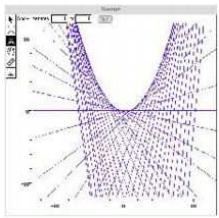

\centerline{\BoxedEPSF{tangents.parabola scaled 400}}
\bigskip
\caption{Each line tangent to the discriminant parabola represents 
the set of quadratic polynomials with a fixed root.} 
\end{figure} 

A simple experiment with a mapping software triggered this investigation.
Figure 1 shows 
the effect of applying the Vi\`{e}te map $\mathcal V$ to a grid of vertical and horizontal lines in the 
$uv$-plane. At the first glance, the result is quite surprising: not only the images 
of lines under \emph{quadratic} map $\mathcal V$ are \emph{lines}, but these lines 
seem to be \emph{tangent} to a parabola! In fact, this parabola $\mathcal D$ 
is the $\mathcal V$-image of the diagonal $\{u = v\}$. Its parametric 
equation is $(b, c) = (-2u, u^2)$. Hence, the equation of $\mathcal D$ is 
quite familiar to the frequent users of the quadratic formula: $b^2 - 4c = 0$. 

In order to understand the tangency phenomenon, consider the 
Jacobi matrix of the Vi\`{e}te map
\[D\mathcal V =  \left( \begin{array}{ccc}
-1 & -1 \\
v & u &
\end{array}
\right).
\] 
Its determinant $J\mathcal V = v - u$. It vanishes along the diagonal line 
$L \subset \A^2_{root}$, where the rank of $D\mathcal V$ drops to 1. This reinforces 
what we already have derived from the symmetry argument: under 
the Vi\`{e}te map, the root plane is ramified over the coefficient plane along 
the discriminant parabola $\mathcal D$.

The kernel of 
$D\mathcal V|_L = Span\{(1,-1)\}$ does not contain the diagonal. The 
$\mathcal V$-image of any smooth curve, which intersects with the diagonal at 
a point $a$ and is transversal there to the kernel of $D\mathcal V|_L$, is tangent to 
the discriminant parabola $\mathcal D = \mathcal V(L)$ at $\mathcal V(a)$. 
In particular, the images of vertical and horizontal lines are tangent to $\mathcal D$. 
However,  $\mathcal V$ maps each vertical line $l_{u_\star} := \{u = u_\star\}$ to a \emph{line} 
$(b, c) = (-u_\star - v,\, u_\star v)$. Therefore, the grid of vertical and horizontal lines 
is mapped by the Vi\`{e}te map to the enveloping family of the discriminant parabola. 

Since the line $\mathcal V(l_{u_\star})$ is the set of all quadratic polynomials with 
a fixed root $u_star$, its points must satisfy the relation $u_\star^2 + bu_\star + c = 0$. 
Therefore, the \emph{slope} of the line 
$\mathcal V(l_{u_\star}) = \{ c = -u_\star b -u_\star^2\}$ is equal 
to \emph{minus the root} $u_\star$! 

As a result, an algebraic problem of solving a quadratic 
equation $z^2 + bz + c = 0$ is equivalent to a geometric problem of finding 
lines passing through the point $(b, c)$ and tangent to the curve $\mathcal D$.
These observations are summarized in 

\begin{prop}
Over the complex numbers, through every point $(b, c) \notin \mathcal D$, 
there are exactly two complex lines tangent to $\mathcal D$. Through every point 
$(b, c) \in \mathcal D$, the tangent line is unique. 

Over the reals, through 
each point of the domain $\mathcal U_+ = \{c < b^2/4\}$, there exists a pair 
of tangent lines, while through each point of the domain $\mathcal U_- = \{c > b^2/4\}$
no such a line exists.

The slopes of these tangents equal to minus the roots of the quadratic equation 
$z^2 + bz + c = 0$.\qed
\end{prop}

\begin{figure}[ht]
\centerline{\BoxedEPSF{SolvingMachine scaled 450}}
\bigskip
\caption{}
\end{figure}

Figure 2 depicts an analog device which is based on this theorem. It solves 
quadratic equations over the field $\R$. The discriminant parabola is modeled by 
a parabolic rim attached to the $bc$-plane. The device consists of two rulers 
hinged by a pin. We solve an equation by placing the pin at the corresponding point
$(b, c)$ and adjusting the rulers to be tangent to the rim. Then we read the 
measurements of their slopes. 
\bigskip

"Completing-the-square" magic calls for a substitution $z \Rightarrow z - t$\; ($t = -b/2$),
 which transforms a given polynomial $P(z) = z^2 + bz + c$ into a polynomial 
$Q(z) = P(z - t)$ of the form  $z^2 - \tilde c$. This kind of substitutions defines 
a $t$-parametric group of transformations 
\begin{eqnarray}
\Phi_t(b, c) = (b - 2t,\; c - bt + t^2) = (P'(-t),\; P(-t))
\end{eqnarray}
in the $bc$-plane. The corresponding transformation in the root plane amounts 
to a simple shift $\Psi_t(u, v) = (u + t,\, v + t)$. In other words, 
$\mathcal V(\Psi_t(u, v)) = \Phi_t(\mathcal V(u, v))$. 

Evidently, $\Psi_t$ preserves the Jacobian $J\mathcal V = v - u$. 
The Jacobian changes sign under the permutation $(u, v) \Rightarrow (v, u)$.
Therefore, it can not be expressed in terms of $b$ and $c$. However, its square 
$(J\mathcal V)^2 = (v - u)^2$ is invariant under the permutation and admits 
a $bc$-formulation as $\Delta(b, c) = b^2 -4c$. Therefore, the \emph{discriminant 
polynomial} $\Delta(b, c)$ must be \emph{invariant} under the $\Phi_t$-flow.  As a result,
the $\Phi_t$-trajectories form the family of parabolas $\{\Delta(b, c) = const\}$. 
In particular, the discriminant parabola is a trajectory of the $\Phi_t$-flow. 

\begin{figure}[ht]
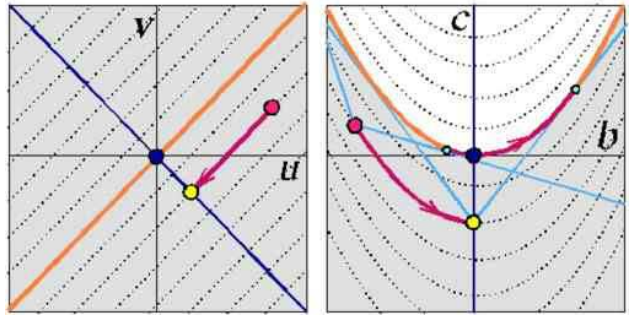

\centerline{\BoxedEPSF{parabolic.flow.1 scaled 450}} 
\bigskip
\caption{The Vi\`{e}te map is equivariant under the flows  $\{\Psi_t\}$ and $\{\Phi_t\}$.  
Transformation $\Phi_t$ acts on the enveloping family of the discriminant parabola by "adding $t$ to 
their slopes".}
\end{figure} 

Each $\Phi_t$-trajectory intersects with the vertical line $\{b = 0\}$ of \emph{reduced} 
quadratic polynomials at a single point. The procedure of completing the square  
amounts to traveling along the trajectory $\Phi_t((a, b))$ until, after $t = -b/2$ 
units of time, it hits the line $\{b = 0\}$ at the point $(0, \, c - b^2/4)$.

We notice that, for a fixed $t$,  $\Phi_t$ is an affine transformation of the 
$bc$-plane. Hence, it maps lines to lines. By the argument above, $\Phi_t$ also 
preserves the family of lines, tangent to the discriminant parabola. 
As the argument suggests, a tangent line with a slope $k$ is mapped by
$\Phi_t$  to  a tangent line with the slope $k + t$.  
In fact, the slopes of lines passing through a point $(0, -d)$ and 
tangent to $\mathcal D$ are $\pm \sqrt d$. 

Therefore, the quadratic formula reflects the following geometric recipe:
\begin{itemize} 
\item apply $\Phi_{-b/2}$ to $(b, c)$ to get to a point $Q = (0,\, c - b^2/4)$.
\item construct the tangent lines to the discriminant parabola through $Q$. 
\item flow them back by $\Phi_{b/2}$ to get tangent lines through $(b, c)$.  
\end{itemize}

Curiously, the flow $\Phi_t$ \emph{preserves the euclidean area} in the $bc$-plane: 
for a fixed $t$,  formula (2.1) describes $\Phi_t$ as a composition of a 
linear transformation $(b, c) \rightarrow (b,\, c - bt)$ with the determinant 1, 
followed by a shift $(b, c) \rightarrow (b - 2t,\, c + t^2)$. 
\bigskip

There is an alternative approach which leads to the same geometric observations
and does not involve the Vi\`{e}te map. However, it  calls for a trip to the 
3rd dimension. 
\begin{figure}[ht]
\centerline{\BoxedEPSF{saddle scaled 600}} 
\bigskip
\caption{}
\end{figure} 
Consider the surface $S = \{z^2 + bz + c = 0\}$ in the $bcz$-space (cf. Figure 4). 
It admits a $bz$-parametrization 
\begin{eqnarray}
\mathcal H: (b, z) \rightarrow (b,\, -bz - z^2,\, z).
\end{eqnarray}
Let $\mathcal F$ be the composition of the parametrization $\mathcal H$ with 
the obvious projection
$\mathcal P : (b, c, z) \rightarrow (b, c).$
It is given by the formula 
\begin{eqnarray}
\mathcal F: (b, z) \rightarrow (b,\, -bz - z^2).
\end{eqnarray}
Evidently, the $z$-function, restricted to the preimage $(\mathcal P|_S)^{-1}((b, c))$, 
gives the roots of $z^2 + bz + c$.  

We focus on the singularities of $\mathcal F$. Its Jacobi matrix is 
\[ D\mathcal F = \left( 
\begin{array}{cc}
1 \, & -z\\
0\, & -b - 2z 
\end{array}
\right).
\]
and its Jacobian $J\mathcal F = -b - 2z$. The rank of $D\mathcal F $ 
drops to 1 when $b  + 2z = 0$, that is, when $P'(z) = 0$. Here $P(z) = 
z^2 + bz + c$. 

Thus, the set of singular points for the projection $\mathcal P|_S$ 
is a curve $C$ in the $bcz$-space given by two equations
\[ 
\begin{cases}
%\left\{ 
\begin{array}{rcc}
z^2 + bz + c &=& 0\\
2z + b &=& 0
\end{array}
\end{cases}
\]
Expelling $z$ from the system we get the equation $c = b^2/4$ of the 
ramification locus for the projection $\mathcal P$ from the surface 
$S$ to the $bc$-plane. Again, as with the Vi\`{e}te map, the discriminant 
parabola is the ramification locus for the projection $\mathcal P$. 
Of course, this is not surprising: the system of equations 
$\{P(z) = 0,\, P'(z) = 0\}$ tells us that the polynomial $P$ has a multiple 
root, that is, $P(z)$ is of the form $(z - u)^2$. In the $bc$-plane, such 
polynomials form the discriminant parabola $\mathcal D$. 
\bigskip

For a fixed a number $u$, let $N^u$ to denote the intersection of the surface $S$ with the plane 
$\{z = u\}$ in the $bcz$-space. This intersection is a \emph{line} 
defined by two equations $\{z = u\}$ and $\{u^2 + bu + c = 0\}$. Hence, $S$ is a \emph{ruled 
surface} comprised of the distinct lines $N^u$. We notice that each line $N^u \subset S$
hits the critical curve $C \subset S$ at a single point:  there is a single 
monic quadratic polynomial with a root $u$ of multiplicity 2.\smallskip 

The projection of $N^u$ in the bc-plane is a line $T^u = \{c + ub + u^2 = 0\}$, 
which, in view of our analysis of the Vi\`{e}te map, is tangent to the discriminant 
parabola. We also can verify this property directly by comparing the 
$\mathcal P$-images of the lines $N^u$ and the line tangent to the $\mathcal P$-critical
curve $C \subset S$ at the point $N^u \cap C$. Thus, the enveloping family of 
the discriminant parabola is the $\mathcal P$-image of the $u$-family of lines 
$\{N^u\}$ comprising $S$.

%%%%%%%%%%%%%%%%%%%%%%%%%%%%%%%%%%%%%%%%%% 

\section{ruled geometry of cubic discriminants}

We build on the observations from Section 2 to investigate the 
discriminant surface for cubic polynomials. This is the simplest case revealing 
the \emph{stratified} ruled nature of the discriminant varieties. 

Facts about the geometry of the cubic discriminant we are going to 
describe here can be found somewhere else (cf. [BG], 5.36). Often we differ 
from these sources only in the interpretation. This 
interpretation will allow us (cf. Sections 6, 7) to investigate the case 
of discriminant varieties for polynomials of any degree.  \smallskip 

Now, our main object of interest is a monic cubic polynomial 
$P(z) = \break z^3 + bz^2 + cz + d$.  
Such polynomials can be coded by points $(b, c, d)$ of the coefficient space  
 $\A^3_{coef}$.\footnote{As always, the coefficient space comes in two flavors: real and complex.}

In order to incorporate the roots of polynomials into the picture, consider the 
hypersurface $S_1$
\begin{eqnarray}
z^3 + bz^2 + cz + d = 0
\end{eqnarray}
in the 4-dimensional space $\A^1\times \A^3_{coef}$ with the cartesian coordinates $z, b, c, d$. 

Put $Q(z, b, c, d) = z^3 + bz^2 + cz + d$. Since the gradient 
$\nabla Q =  (P'(z),\, z^2,\, z,\, 1) \neq 0$,  the hypersurface $S_1$ 
is non-singular. It can be viewed as the graph of the function 
$d(z, b, c) = - z^3 - bz^2 - cz$. Therefore, $S_1$ admits a $(z, b, c)$-parametrization 
by a 1-to-1 polynomial map  
\begin{eqnarray}
\mathcal H_1: (z, b, c) \rightarrow (z, b, c, d) = (z,\, b,\, c, - z^3 - bz^2 - cz).
\end{eqnarray}

Denote by $\mathcal P$ the projection $(z, b, c, d) \rightarrow (b, c, d).$ 
Our immediate goal is to analyze the singularities of this projection, being 
restricted to the hypersurface $S_1$.  That is, we will investigate the singularities 
of the composition $\mathcal F_1 = \mathcal P \circ \mathcal H_1$ given by 
\begin{eqnarray}
\mathcal F_1: (z, b, c) \rightarrow (b, c, d) = (b,\, c, - z^3 - bz^2 - cz).
\end{eqnarray}

The Jacobi matrix $D\mathcal F_1$ of $\mathcal F_1$ is of the form
\[\left( 
\begin{array}{ccc}
0 \, & 0 \, & -3z^2 -2bz -c \\
1 \, & 0 \, & -z^2\\
0\, & 1 \, & -z
\end{array}
\right).
\]  
Unless the derivative $P'(z) = 3z^2 + 2bz + c = 0$, the rank of $D\mathcal F_1$ 
is 3. When $P'(z) = 0$, it drops to 2. Thus, a point $(z, b, c, d) \in S_1$ is 
singular for the projection $\mathcal P|_{S_1}$, if and only if, two conditions are satisfied: 
$P(z) = 0$ and $P'(z) = 0$.  This happens exactly when $z$ is a root of $P$ of 
multiplicity $\geq 2$. Therefore, the singular locus $S_2$ of $\mathcal P|_{S_1}$  in 
the $zbcd$-space is the intersection of two hypersurfaces
\begin{eqnarray}
z^3 + bz^2 + cz + d = 0\nonumber\\
3z^2 + 2bz + c = 0.
\end{eqnarray}
Solving this system for $c$ and $d$, the non-singular surface $S_2$ can be parametrized 
by $z$ and $b$:
\begin{eqnarray}
\mathcal H_2: (z, b) \rightarrow (z, b, c, d) = (z,\, b, -3z^2 -2bz,\, 2z^3 + bz^2).
\end{eqnarray}
Composing $\mathcal H_2$ with the projection, we get: 
\begin{eqnarray}
\mathcal F_2: (z, b) \rightarrow (b, c, d) = (b, -3z^2 -2bz,\, 2z^3 + bz^2).
\end{eqnarray}

The Jacobi matrix $D\mathcal F_2$ of $\mathcal F_2$ is
\[\left( 
\begin{array}{ccc}
0 \, & -6z -2b \, & 6z^2 + 2bz \\
1 \, & -2z \, & z^2
\end{array}
\right)
\] 
Generically, it is of rank is 2. The rank drops to 1 when $P''(z) = 6z + 2b = 0$, 
that is, when $(z, b, c, d)$ belongs to a curve $S_3 \subset S_2$, defined by the three 
equations 
\begin{eqnarray}
P(z) &=& z^3 + bz^2 + cz + d = 0\nonumber\\
P'(z) &=& 3z^2 + 2bz + c = 0\nonumber\\
P''(z) &=& 6z + 2b = 0
\end{eqnarray}
The curve $S_3$ admits a parametrization by $z$
\begin{eqnarray}
\mathcal H_3: z \rightarrow (z, b, c, d) = (z, -3z,\, 3z^2, -z^3).
\end{eqnarray}
Its non-singular projection $\mathcal F_3$ into the $bcd$-space is given by 
\begin{eqnarray}
\mathcal F_3: z \rightarrow (b, c, d) = (-3z,\, 3z^2, -z^3),
\end{eqnarray}
In view of (3.7), this curve $\mathcal D_3$ represents cubic polynomials with a single root of 
multiplicity 3. 
\smallskip

Let $\mathcal D_2$ denote the image $\mathcal P(S_2)$ of the surface $S_2$, and 
$\mathcal D_3$---the image $\mathcal P(S_3)$ of the curve $S_3$ in the $bcd$-space.
For the reasons, that will be even more apparent in Section 4, we call these 
images the \emph{discriminant surface} and the \emph{discriminant curve}. 
By the definitions, $\mathcal D_3 \subset \mathcal D_2 \subset \A^3_{coef}$. 
Figures 5 and 6 show the discriminant surface from different points of view.
\bigskip 

\begin{figure}[ht]
\centerline{\BoxedEPSF{D^2 scaled 550}}
\bigskip
\caption{$\mathcal D_2$ is a ruled surface comprised of lines tangent to 
the discriminant curve $\mathcal D_3$ (perceived as a loop).}
\end{figure} 

\begin{figure}[ht]
\centerline{\BoxedEPSF{D^2.1 scaled 800}} 
\bigskip
\caption{This view of the discriminant surface $\mathcal D_2$ reveals its symmetry with 
respect to the involution $(b, c, d) \rightarrow (-b,\, c, -d)$.}
\end{figure}

Looking from a point $P = (b, c, d) \in \A^3_{coef}$ against the projection $\mathcal P$, 
we see all the points of the hypersurface $S_1$ (cf. (3.1)) suspended over $P$,
in other words, all the roots of the polynomial $P(z) = z^3 + bz^2 + cz + d$. 
Therefore, 
the preimage $S_1 \cap \mathcal P^{-1}(A)$ can contain 1, 2, or 3 points.
Over the complex numbers, the cardinality of the preimage equals 3 when 
$P \in \C^3 \setminus \mathcal D_2$, it is 2 when $P \in \mathcal D_2 \setminus \mathcal D_3$, 
and 1 when 
$P \in \mathcal D_3$. Thus, $\mathcal P|_{S_1}$ is 3-to-1 
map, \emph{ramified} over the discriminant surface $\mathcal D_2$. Similarly, 
$\mathcal P|_{S_2}$ is 2-to-1 map, ramified over the discriminant curve $\mathcal  D_3$. 
Finally, $\mathcal P|_{S_3}$ is 1-to-1 map. 

Over the real numbers, the situation is more complex: the 
surface $\mathcal D_2$ divides $\R^3_{coef}$ into chambers, and, by the implicit 
function theorem, the 
cardinality of $S_1 \cap \mathcal P^{-1}(A)$ remains constant for all the $P$'s 
in the interior of each chamber. Since a real cubic polynomial with no 
multiple roots has, alternatively, three real roots or a single real root, 
the chambers can be only of two types. (In the next section we will see that 
actually, there is a single chamber of each type.) 
We notice that, if a \emph{real} 
cubic polynomial has a multiple complex root, then all its roots are real.
Therefore, even over the reals, for $P \in \mathcal D_2 \setminus \mathcal D_3$, the 
preimage $S_1 \cap \mathcal P^{-1}(A)$ consists of two points. 
\bigskip

Now, let's return to the $zbcd$-space. 
For a given number $u$, consider the intersection $N^u_1$ of the hyperplane 
$\{z = u\}$ with the hypersurface $S_1$. This intersection selects 
all the quadruples $(u, b, c, d)$ with the property 
$P(u) = u^3 + bu^2 + cu + d = 0$, where $u$ is  fixed. 
 The $\mathcal P$-image of $N^u_1$ is the 
surface $T^u_1$ of  all cubic polynomials with the number $u$ as their 
common root. Evidently, the map $\mathcal P: N^u_1 \rightarrow T^u_1$ is 1-to-1 and 
onto. 
One thing is instantly clear: $u^3 + bu^2 + cu + d = 0$ defines 
a \emph{linear} relation among $b, c, d$---an \emph{affine plane} in $\A^3_{coef}$.
Scanning by $u$, we see that the hypersurface $S_1$ is a disjoint 
union of its $u$-slices --- the planes $N^u_1$, 
i.e. it is a \emph{ruled} hypersurface.

What is also clear, that each plane $N^u_1$ hits the surface $S_2$ (of 
polynomials with multiple roots) along a \emph{line} $N^u_2$, defined by the
equations (3.4) and the equation $\{z = u\}$. Indeed, the set of cubic equations 
with a root $u$ contains the  set of cubic equations with a root $u$ of 
multiplicity $\geq 2$. In turn, the line $N^u_1$ hits the curve $S_3$ at a single 
point: there is a single monic cubic polynomial with the root $u$ of multiplicity 
3. Therefore, the surface $S_2$ is a ruled surface comprised of disjoint lines
$N^u_2$. Since $\mathcal D_2 = \mathcal P(S_2)$ and $\mathcal P$ maps lines to 
lines, $\mathcal D_2$ is also a ruled surface comprised of lines defined by 
\begin{eqnarray}
u^2b + uc + d& = - u^3\nonumber\\
2ub + c& = - 3u^2.
\end{eqnarray}

\smallskip 

Let's concentrate on the case when $u$ is a root of multiplicity $\geq 2$.
Consider a plane $N^u_1$ through a point 
$(u, P)$, $P = (b, c, d)$, of the surface $S_2 \subset S_1$, and the plane $\tau_{(u,P)}$ 
\emph{tangent} to $S_2$ at $(u, P)$.  We will show that $T^u_1 = \mathcal P(N^u_1)$ 
is \emph{tangent} to $\mathcal D_2$  at the point $P$. It will suffice 
to check that vectors, tangent to $S_2$ at $(u, P)$ project into the plane  $T^u_1$. 

Using (3.4), the plane $\nu_{(z, P)}$, normal to $S_2$ at a point $(z, P)$, is spanned 
by two gradient vectors $\nabla_1(z, P) = (P'(z), z^2, z, 1)$ and 
$\nabla_2(z, P) = (P''(z), 2z, 1, 0)$. Note that, when $(z, P) \in S_2$, then 
$\nabla_1(z, P) = (0, z^2, z, 1)$. Denote by $n = n(z)$ the vector $(z^2, z, 1)$. 
In the new notation,  $\nabla_1(z, P) = (P'(z), n(z))$ and $\nabla_2(z, P) = (P''(z), n'(z))$.

Any vector $(a, v) \in \A^1\times\A^3$, tangent to $S_2$ at $(u, P)$  must be orthogonal 
to $\nabla_1(u, P) = (0,\, n(u))$ and $\nabla_2(u, P) = (P''(u),\, n'(u))$, in other words, 
$(a, v)$ must satisfy the system
\begin{eqnarray}
v \bullet n(u)& =& 0\nonumber \\
v \bullet n'(u)& =& - a P''(u),
\end{eqnarray}
where "$\bullet$" stands for the scalar product. We notice that if $v$ satisfies the 
first equation in (3.11), then one can always to find an appropriate $a$, provided 
$P''(u) \neq 0$. On the other hand, when $P''(u) = 0$, that is, when $(u, P) \in S_3$,
then (3.11) collapses to $\{v \bullet n(u) = 0, \; v \bullet n'(u) = 0\}$ and the $a$ 
is free. In such a case, $v$ must belong to a \emph{line} $L^u$ passing through $P$. 
This line is the $\mathcal P$-image of the tangent plane $\tau_{(u,P)}$ and 
its parametric equation is of the form $\{P + v\}_v$, where $v$ is a subject to 
the orthogonality conditions $\{v \bullet n(u) = 0, \; v \bullet n'(u) = 0\}$.   

Therefore, if $P \in \mathcal D_2^\circ := \mathcal D_2 \setminus \mathcal
D_3$, then any vector $v$ (with its origin at $P$) orthogonal to $n(u)$ belongs to the plane
$\mathcal P(\tau_{(u,P)})$. As a result, such a $v$ must be tangent to 
$\mathcal D_2 = \mathcal P(S_2)$ at $P$. On the other hand, the vector $n(u)$ is 
normal to the plane  $T^u_1 := \{u^2b + uc + d = -u^3\}$ passing through $P$. 
So, as  affine planes, $T^u_1 = \mathcal P(\tau_{(u,P)})$, and therefore, $T^u_1$ 
must be tangent to $\mathcal D_2^\circ$ at $P$, provided $P(u) = 0$. 

Note that, for any $P \in \mathcal D_2$, there is a single point $(u, P)$
in $\mathcal P^{-1}(P) \cap S_2$: a cubic polynomial can not have more 
than one multiple root. Therefore, $\mathcal P: S_2^\circ \rightarrow 
\mathcal D_2^\circ$ is a regular embedding.
\smallskip

In the same spirit, one can check that when $(u, P) \in S_3$, the line 
$\{P + v\}_v = \mathcal P(\tau_{(u,P)})$ determined by 
$\{v \bullet n(u) = 0, \; v \bullet n'(u) = 0\}$ coincides the line 
$T^u_2 \subset T^u_1$, defined by two equations 
$\{P(u) = 0,\, P'(u) = 0\}$ as in (3.10). By its 
definition, $T^u_2 \subset \mathcal D_2$. Furthermore, each line $T^u_2$ 
is \emph{tangent} to the curve $\mathcal D_3$ at their intersection point $P^u$ which  
corresponds to a polynomial of the form $P(z) = (z - u)^3$. Indeed, the 
vector $w(u) = (-3,\, 6u, -3u^2)$, tangent to the curve $\mathcal D_3$ at $P^u$, 
is orthogonal to $n(u)$ and $n'(u)$. This becomes evident using the identities 
$\partial_u\{(z - u)^3\} = - 3z^2 + 6zu^2 - 3u^2 = (z^2, z, 1)\bullet(-3,\, 6u, -3u^2) = n(z)\bullet w(u)$  
and $\partial_z\partial_u\{(z - u)^3\} = - 6z + 6u^2  = (2z, 1, 0)\bullet(-3,\, 6u, -3u^2) = n'(z)\bullet w(u)$
--- just substitute $z = u$.\smallskip

One can check that, for \emph{distinct} $u$, the systems (3.10) do not share a common 
solution $(b, c, d)$, in other words, all the lines $T^u_2$ are disjoint.  In combination 
with the previous arguments this leads to a conclusion which could be predicted examining the 
images in Figures 5 and 6.

\begin{prop} The discriminant surface $\mathcal D_2$ is a ruled surface comprised 
of the disjoint lines $T^u_2 \subset \A^3_{coef}$ defined by two constraints $P(u) = 0,\, P'(u) = 0$ 
as in (3.10). Each line $T^u_2$ is  tangent to the discriminant 
curve $\mathcal D_3$ at a point $P^u$, corresponding to the polynomial $P(z) = (z - u)^3$. 
In other words, $\mathcal D_2$ is spanned by lines tangent to $\mathcal D_3$. \qed
\end{prop}

The same conclusion can be reached following a different approach. Both treatments 
will play different and complementary roles in Section 6. 

The curve $\mathcal D_3$ 
admits a parametrization $A(u) = (-3u,\, 3u^2, -u^3)$.  
The velocity vector $w(u)$  at $A(u)$ is equal to 
$\dot{A}(u) = (-3,\, 6u, -3u^2)$. A line, tangent to  $\mathcal D_3$ at $A(u)$, has a 
$t$-parametric equation $A(u) + t \dot{A}(u)$. This line corresponds to the $t$-family 
of monic cubic polynomials of the form $P(z) - t P'(z) = (z -u)^3 - 3t(z - u)^2$. 

A generic polynomial $Q(z) = (z - u)^2(z - u')$ in $\mathcal D_2$, for 
an appropriate choice of $t$, can be represented in the form $P(z) + t P'(z)$ . 
To do it, we need to solve for $t$ the $z$-functional equation  
$(z - u)^2(z - u') = (z - u)^3 + 3t(z - u)^2$. Miraculously, it has a unique solution 
$t = (u' - u)/3$ ! 
Therefore, any point $Q$ in $\mathcal D_2$ lies on a line $l^Q$ tangent to $\mathcal D_3$. 
At the same time, an attempt 
to solve for $t$ the $z$-functional equation $(z - u)^2(z - u') = (z - u'')^3 + 3t(z - u'')^2$ fails 
when $u'' \neq u$.  Thus, the tangent line to $\mathcal D_3$ through $Q$ is unique.\qed
\bigskip

Let us make a few crucial observations about 
the way in which a plane can be tangent to the ruled surface $\mathcal D_2$. 
First we notice that, if a ruled surface $\mathcal D$ has a tangent plane $T$ 
at a non-singular point $P \in \mathcal D$, then $T$ must contain all the lines 
through $P$ from the family which forms $\mathcal D$. Therefore,  
any plane $T$, tangent to $\mathcal  D_2$ at $P \in \mathcal  D_2^\circ$, contains 
the unique line $T^u_2$ through $P$. In turn, $T^u_2$ is tangent to $\mathcal D_3$ at a different 
point $P^u$. This geometry is depicted in Figure 7.\smallskip 

\begin{figure}[ht]
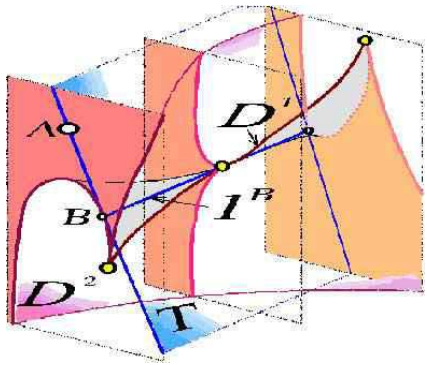

\centerline{\BoxedEPSF{tangent.plane scaled 450}} 
\caption{The slices of $ \mathcal D^2$ by planes $\{b = const\}$. 
Any plane $T$ through a point
$A$ and  tangent to $\mathcal D_2$ is tangent along the whole line $l^B$. In turn, $l^B$ is 
tangent to the curve $\mathcal D_3$. }
\end{figure}

We have seen that $\mathcal D_2^\circ$ is a smooth surface. However, 
the surface $\mathcal D_2$ fails to be smooth at the points of the discriminant 
curve $\mathcal D_3$: it has a cusp-shaped fold along $\mathcal D_3$ (cf. Corollary 5.1 and Figure 8).  
Therefore, we need to clarify the notion of a "tangent" plane to  $\mathcal  D_2$ at 
the points of its singular locus $\mathcal  D_3$.\smallskip 

We define the  "tangent" plane to $\mathcal D_2$ at $P \in \mathcal D_3$ to be  
the unique plane spanned  by the velocity and acceleration vectors at $P$ of the parametric curve 
$\mathcal D_3$---the, so called,  \emph{osculating} plane of the curve. The osculating 
plane at $P \in \mathcal D_3$ happens to be the \emph{limit}, as $Q$ approaches $P$, 
of tangent planes at smooth points $Q \in \mathcal  D_2^\circ$.  
This is another small miracle of the discriminant surface: although it is singular 
along $\mathcal D_3$, the tangent planes of smooth points stabilize towards $\mathcal D_3$ 
(cf. Proposition 3.2)---the tangent bundle of $\mathcal  D_2^\circ$ extends to a bundle over 
$\mathcal  D_2$. 
\smallskip

In order to verify these claims, consider a $ts$-parametric equation of an osculating plane 
through a point $A(u) = (-3u,\, 3u^2, -u^3)$ on the discriminant curve $\mathcal D_3$ 
--- a plane which is spanned by the velocity $\dot A(u)$ and the acceleration $\ddot A(u)$  
vectors:
\begin{eqnarray}
(b, c, d) = (-3u,\, 3u^2, -u^3) + t(-3,\, 6u, -3u^2) + s(0,\, 6, -6u)  
\end{eqnarray}
Expelling $t$ and $s$ from these three equations, we get a somewhat familiar 
relation between $b, c, d$ and $u$:\, $\{u^3 + bu^2 + cu + d = 0\}$! Conversely, 
if $u$ is a root of an equation $z^3 + bz^2 + cz + d = 0$, then $(b, c, d)$ 
belongs to the osculating plane in (3.12) at $A(u)$: just put $t = u + \frac{1}{3}b$ and 
$s = \frac{1}{6}c  + \frac{1}{3}bu + \frac{1}{2}u^2$.  

We notice that the vector
$n(u) = (u^2,\, u,\, 1)$, normal to $\mathcal  D_2^\circ$ at the points of the 
line $T^u_2$ , is also normal to both vectors: 
$\dot A(u) = (-3,\, 6u, -3u^2)$, $\ddot A(u) = (0,\, 6, -6u)$. 
Therefore, the osculating plane at $A(u) \in \mathcal D_3$ coincides with 
the affine plane $T^u_1$ tangent to $\mathcal D_2^\circ$ along the line $T^u_2$.
We have proved the following
 
\begin{prop} If $u$ is a root of the polynomial $P(z) = z^3 + bz^2 + cz + d$, then the 
osculating plane of the curve $\mathcal  D_3$ at the point $A(u)  = (-3u,\, 3u^2, -u^3)$ 
contains the point $(b, c, d)$. That plane coincides with the affine plane 
$T^u_1$ and thus, is tangent to the surface $\mathcal  D_2^\circ$ 
along the line $T^u_2$. In turn, $T^u_2$ is tangent to the curve at $A(u)$. \qed  
\end{prop}

The embedding $\mathcal D_2 \subset \A^3_{coef}$ has a characteristic property 
described in
\begin{prop}
\emph{Any} plane in $\A^3_{coef}$, passing through the point $P = (b, c, d)$ and tangent 
to the surface $\mathcal D_2$, is of the form $T^u_1 := \{bu^2 + cu + d = -u^3\}$, 
where $u$ is a root of the polynomial $P(z) = z^3 + bz^2 + cz + d$
\end{prop}

\textbf{Proof.}\quad  The considerations above already contain 
the proof. We have seen that the planes $\{T^u_1\}$ are exactly the 
planes tangent to $\mathcal  D_2^\circ$. Moreover, by Proposition 3.2, the tangent cones of 
$\mathcal D_2$ at the points of $\mathcal D_3$ belong to the same family of 
planes. On the other hand, every point $P \in \A^3_{coef}$ belongs to 
each of the planes $T^u_1$, where $u$ ranges over the roots of $P(z)$.\qed
\smallskip

Even the existence of \emph{finitely many} tangent planes through a generic $P$ is an
extraordinary fact: for a \emph{general} surface $S$, there exists 
an 1-parametric family of planes passing through $P$ and tangent to $S$. 
What distinguishes $\mathcal D_2$ from a general surface, is its ruled 
geometry--- $\mathcal D_2$ is formed by the lines tangent
to the spatial curve $\mathcal D_3$. Via the Gaussian map, the tangent planes of the surface 
$\mathcal D_2$ form an 1-dimensional set in the Grassmanian $Gr(3,2) \approx \P^2$, 
while a generic surface generates a 2-dimensional Gaussian image.
\bigskip

Now we are ready to restate a fundamental relation between the ruled stratified 
geometry of the determinant variety $\mathcal D_2$ and the roots of the cubic 
monic polynomials. The proposition below is 
just a repackaging of the propositions that already have been established. 

We start with the complex case which, as usual, is more uniform. 
\begin{thm} 

\begin{itemize} 
\item  
Through any point $P$ of the stratum $\C^3_{coef} \setminus \mathcal D_2$, 
there are exactly 3 planes, tangent to the discriminant surface $\mathcal D_2$. 

\item  
Through any point $P$ of the stratum $\mathcal D_2^\circ$, 
there are exactly 2 planes, tangent to  
the surface $\mathcal  D_2$. One of 
these planes is tangent to $\mathcal D_2$ along a line passing 
through $P$.

\item Finally, through any point $P$ of $\mathcal  D_3$, 
there is a single plane, tangent to 
the surface $\mathcal  D_2$. It 
is the osculating plane of the curve  $\mathcal  D_3$ at $P$.
\end{itemize}

Each of these planes $T^u_1$ is  tangent to the surface 
$\mathcal  D_2^\circ$ along a line $T^u_2 \subset \mathcal  D_2$. 
In turn,  the line is tangent to the discriminant curve $\mathcal  D_3$.
\smallskip

Moreover, if $P = (b, c, d) \in \C^3_{coef}$, then each 
of the tangent planes $T^u_1$ passing through $P$ is described  
by an equation of the form $\{u^2 b + uc + d  = -u^3\}$, where  $u$ 
runs over the distinct complex roots of the polynomial 
$P(z) = z^3 + bz^2 + cz + d$. \smallskip

In fact, each affine plane $T^u_1$ is the osculating plane of 
the discriminant curve at the point of $\mathcal  D_3$ corresponding 
to the polynomial $P(z) = (z - u)^3$. \qed
\end{thm}

The case of cubic polynomials with real coefficients is similar, but has 
a bit more structure and complexity. At the same time, 
the proof is virtually the same, as in the complex case. The stratum 
$\R^3_{coef} \setminus \mathcal D_2$ is divided into two chambers $\mathcal U_3$ and 
$\mathcal U_1$: the first corresponds to real cubic polynomials with 3 distinct 
real roots, the second---with a single simple real root. 
\begin{thm} 
\begin{itemize} 
\item  
Through any point $P \in \mathcal U_3$, 
there are exactly 3 planes, tangent to the discriminant surface $\mathcal D_2$. 
\item 
Through any point $P \in \mathcal U_1$, 
there is exactly 1 plane, tangent to  $\mathcal D_2$.
\item  
Through any point $P \in \mathcal D_2^\circ$, 
there are exactly 2 planes, tangent to  
 $\mathcal  D_2$. One of 
these planes is tangent to $\mathcal D_2$ along a line passing 
through $P$.

\item Finally, through any point $P$ of $\mathcal  D_3$, 
there is a single plane, tangent to 
the surface $\mathcal  D_2$. It 
is the osculating plane of the curve  $\mathcal  D_3$ at $P$.
\end{itemize}

Each of these planes $T^u_1$ is  tangent to the surface 
$\mathcal  D_2^\circ$ along a line $T^u_2 \subset \mathcal  D_2$. 
In turn,  the line is tangent to the discriminant curve $\mathcal  D_3$.
\smallskip

Moreover, if $P = (b, c, d) \in \R^3_{coef}$, then each 
of the tangent planes $T^u_1$ passing through $P$ is described  
by an equation of the form $\{u^2 b + uc + d  = -u^3\}$, where  $u$ 
runs over the distinct real roots of the polynomial 
$P(z) = z^3 + bz^2 + cz + d$. \smallskip

In fact, each affine plane $T^u_1$ is the osculating plane of 
the discriminant curve at the point of $\mathcal  D_3$ corresponding 
to the polynomial $P(z) = (z - u)^3$. \qed
\end{thm}

\begin{cor} \textbf{(Cardano's formula "via tangents")}

It is possible to reconstruct all 
the roots of an equation  $z^3 + bz^2 + cz + d = 0$ 
from the  planes passing through the point $P  = (b, c, d) \in \A^3_{coef}$ 
and tangent to the stratified discriminant pair $\mathcal  D_2 \supset \mathcal  D_3$
--- the tangent planes trough $P$ "solve" the cubic equation. 
Specifically, pick a vector normal to such a tangent plane and 
having the $d$-coordinate 1. Then its $c$-coordinate delivers the 
corresponding root.

In particular, the tangent planes through $P = (0, 0, d)$ have normal
vectors $(\xi^2, \xi, 1)$, where $\xi = \sqrt[3]d$ (complex or 
real).
\qed
\end{cor}

Let's take a flight over the discriminant surface to admire its 
triangular horizon. First, we need a few definitions to 
inform the trip.\smallskip

Given a smooth surface $\mathcal S$ in $\C^3$ and a point $x$ outside $\mathcal S$, one can
associate to each point $y \in \mathcal S$ the unique line $l_{x, y}$ through $x$ and $y$. 
This  defines a map $\pi_x : \mathcal S \rightarrow \P^2_x$ into the projective space $\P^2_x$ 
of lines through $x$. 
We consider the (Zariski) closure $hor(\mathcal S, x)$ of the set of critical points 
for the projection $\pi_x$ 
and call it \emph{the horizon} of $\mathcal S$ at $x$.  Its interior is formed by points 
$y \in \mathcal S$ for which the line $l_{x, y}$ is tangent to $\mathcal S$ at $y$. 
The $\pi_x$-image of $hor(\mathcal S, x)$ in $\subset \P^2_x$, denoted $hor_{\pi}(\mathcal S, x)$,  
is called \emph{the projective horizon} of $\mathcal S$ at $x$. 
Over the real numbers, one gets a refined version of these constructions 
and notions by replacing the space of lines $\P^2_x$ through $x$ by the space of 
rays. This has an effect of replacing the projective plane by the sphere $S^2_x$.

If a surface $\mathcal S$ has a singular locus $K$, then we define $hor(\mathcal S, x)$ and 
$hor_{\pi}(\mathcal S, x)$ as the (Zariski) closure of $hor(\mathcal S \setminus K, x)$ and 
$hor_{\pi}(\mathcal S \setminus K, x)$ in $\mathcal S$ and $\P^2_x$ respectively.
\smallskip

Recall, that the discriminant surface $\mathcal D_2$ has 
a distinct property: if a plane $T = T^u_1$ is tangent to it at a point $Q$, then it is 
tangent to the surface along the entire line $l^Q = T^u_2$ which passes through $Q$. 
Therefore, if $Q \in  hor(\mathcal D_2, P)$, then the line $l^Q \subset  hor(\mathcal D_2, P)$ 
and the projective line $\pi_Q(l^Q) \subset  hor_{\pi}(\mathcal D_2, P)$.
Over the reals, $\pi_Q(l^Q)$ it is a big circle.

We notice that the "naked singularity" $\mathcal D_3$, visible from any  
point $P$ in the coefficient space, it \emph{tangent} to the perceived 
singularity---the projective horizon.
In the complex case if $P \notin \mathcal D^2$, or in the real case when 
$P \in \mathcal U_3$, the discriminant 
curve $\mathcal D_3$ is tangent to the horizon at three points (belonging to 
the three distinct lines which form the horizon). In the real case, when $P \in \mathcal U_1$, 
the discriminant curve is tangent to the horizon line at a single point.
Thus, over the complex numbers, the plane curve $\pi_x(\mathcal D_3) \subset \P^2_x$ is 
\emph{inscribed} in the 
triangular projective horizon. A similar property holds in the real case when 
$P \in \mathcal U_3$. Hence, another distinct property of the discriminant surface:

\begin{cor} For any point $P \in \C^3_{coef} \setminus \mathcal D^2$, the horizon 
$hor(\mathcal D_2, P)$ consists of three lines in a general position in $\C^3_{coef}$.  
The projective horizon $hor_{\pi}(\mathcal D_2, P)$ 
is a union of three projective lines occupying a general position in $\P^2_P$.
The spatial curve $\mathcal D_3$ is inscribed in $hor(\mathcal D_2, P)$, 
while the plane curve $\pi_P(\mathcal D_3)$ is inscribed in the
"triangular" projective  horizon $hor_{\pi}(\mathcal D_2, P)$. 
\smallskip

For any point $P \in \mathcal U_3$, the horizon 
$hor(\mathcal D^2, P)$ consists of three lines in a general position in $\R^3$.  
The projective horizon $hor_{\pi}(\mathcal D^2, P)$ 
is a union of three big circles  occupying a general position in $S^2_P$.

For any point $P \in \mathcal U_1$, the horizon 
$hor(\mathcal D^2, P)$ consists of a single line in $\R^3$,  
while the projective horizon $hor_{\pi}(\mathcal D^2, P)$ 
is a  big circle in $S^2_x$.

The spatial curve $\mathcal D_3$ is inscribed in $hor(\mathcal D^2, P)$, 
while the plane curve $\pi_P(\mathcal D_3)$ is inscribed in $hor_{\pi}(\mathcal D_2, P)$.
\qed
\end{cor}

In general, one might conjecture that the degree of a spatial algebraic curve $\mathcal C$
is perceived as the number of lines in the generic projective horizon of the surface spanned 
by the lines tangent to $\mathcal C$. 
 
%%%%%%%%%%%%%%%%%%%%%%%%%%%%%

\section{The cubic Vi\`{e}te map}

All the results of Section 3 can be understood from a different perspective. 
In Section 2, we described the geometry of the quadratic Vi\`{e}te map. Now we will 
investigate the geometry of the cubic Vi\`{e}te Map.
\smallskip

Let $u, v, w$ be complex roots of a cubic polynomial $P(z) = z^3 + bz^2 + cz + d$. 
Then $P(z) = (z - u)(z - v)(z - w)$. Multiplying the three linear terms, we get 
$P(z) = z^3 -(u + v + w)z^2 + (uv + vw + wu)z - uvw$.  This gives the Vi\`{e}te 
formulas 
\begin{eqnarray}
b& = & -u - v - w  \nonumber\\
c& = & uv + vw + wu\nonumber\\
d& = & -uvw,
\end{eqnarray}
linking roots to coefficients. We think about (3.1) as giving rise to 
a \emph{polynomial map} $\mathcal V$ from the $uvw$-root space $\A^3_{root}$ 
to the $bcd$-coefficient space $\A^3_{coef}$. We call it the \emph{Vi\`{e}te map}.
\bigskip 

By the Fundamental Theorem of Algebra, for any triple $(b, c, d)$ there exists 
a triple of complex numbers $(u, v, w)$, which satisfies the system (4.1), 
in other words, the \emph{complex} Vi\`{e}te map is \emph{onto}. This is not the case for 
the real Vi\`{e}te map.

Because 
the  factorization of $P(z)$ into a product of monic linear polynomials is 
unique up to their ordering, the triple $(b, c, d)$ determines the triple 
$(u, v, w)$ up to permutations in three letters. They form a permutation 
group $S_3$ of order six. Generically, over the complex numbers, the preimage 
$\mathcal V^{-1}(b, c, d)$ 
consists of 6 elements. This happens when $(b, c, d) = \mathcal V(u, v, w)$ with 
$u, v, w$ being distinct. When two of the roots coincide (that is, when the roots 
of $P(z)$ are of multiplicities 1 and 2), $\mathcal V^{-1}(b, c, d)$ consists 
of three elements. Finally, when the polynomial has a single root of multiplicity 
3, $\mathcal V^{-1}(b, c, d)$ is a singleton.
\bigskip

The Fundamental Theorem of Algebra has a fancy formulation in terms of symmetric 
products of the space $\C$ (or even better, of the projective space $\P_1$).
\smallskip

Recall, that the $n$-th symmetric product $S^nX$ of a set $X$ is defined to be 
the $n$-th cartesian product $X^n$ of $X$, divided by the natural action of the 
symmetry group $S_n$. In other words, while points of $X^n$ are \emph{ordered} 
$n$-tuples of points from $X$, points of $S^nX$ are \emph{unordered} $n$-tuples.

In these terms, the algebraic root-to-coefficient map $\mathcal V$  establishes 
an 1-to-1 and onto correspondence 
$\tilde{\mathcal V} : S^n\C_{root} \rightarrow \C^n_{coef}$. In particular,
via the Vi\`{e}te map  $\mathcal V$ in (4.1),  $S^3\C$ and $\C^3$ are isomorphic sets.

The obvious forgetful map $f: \C^n \rightarrow S^n\C$, which strips an ordered $n$-tuple 
of its order, generically, is ($n$!)-to-1. 
\bigskip  

The embedding  $\mathcal D_3 \subset \A^3_{coef}$ provides us with a very 
geometric way of interpreting the Vi\`{e}te map. This interpretation is based on 
Proposition 3.2 and Theorems 3.1, 3.2.  

The discriminant curve $\mathcal D_3$ is a rational curve. It admits  a 1-to-1
parametrization $A = \mathcal F_3 : \A^1 \rightarrow \mathcal D_3$ as in (3.9). 
Given any \emph{unordered} triple of distinct points $A(u), A(v), A(w) \in \mathcal D_3$, 
the corresponding osculating planes $T^u_1,T^v_1, T^w_1$ of the curve at $A(u), A(v), A(w)$ 
all intersect at  
a singleton $\Psi(A(u), A(v), A(w))$ representing the polynomial $P(z) = (z - u)(z - v)(z - w)$.
If $u = w$, we define $\Psi(A(u), A(v), A(u))$ to be the singleton where the osculating plane
$T^v_1$ hits the tangent line $T^u_2$. Of course, this point on $\mathcal D_2^\circ$ corresponds
to the polynomial $P(z) = (z - u)^2(z - v)$. Finally, when $u = v = w$, $\Psi(A(u), A(u), A(u))$   
 is defined to be $A(u)$ which corresponds to $P(z) = (z - u)^3$.
 
This gives rise to well-defined algebraic map $\Psi : S^3\mathcal D_3 \rightarrow \C^3_{coef}$ 
from the symmetric cube of the discriminant curve \footnote{Points of $S^3(\mathcal D_3)$ are 
effective divisors of degree 3 on the curve $\mathcal D_3$.} onto the coefficient space.   

\begin{thm} The Vi\`{e}te map $\mathcal V : \C^3_{root} \rightarrow \C^3_{coef}$ is a composition
of  the forgetful map $f: \C^3_{root} \rightarrow S^3\C_{root}$, the obvious $A$-parametrization map
$S^3A : S^3\C_{root} \rightarrow S^3\mathcal D_3$, 
and the osculating planes map $\Psi : S^3\mathcal D_3 \rightarrow \C^3_{coef}$. 
All the three maps are onto and the maps $S^3A$ and $\Psi$ are 1-to-1. \qed
\end{thm}

In order to describe a crude geometry of the Vi\`{e}te map, we shall concentrate 
on the loci in the coefficient space, where the cardinality $|\mathcal V^{-1}(b, c, d)|$ 
of the preimage   $\mathcal V^{-1}(b, c, d)$ jumps, that is, on the \emph{ramification} loci. 
The previous argument tells us that, over the complex numbers, the condition 
$|\mathcal V^{-1}(b, c, d)| = 1$ 
picks the set of polynomials with a single root of multiplicity 3, 
the condition $|\mathcal V^{-1}(b, c, d)| = 3$ 
picks the set of polynomials with one root of multiplicity 2, finally, 
the condition $|\mathcal V^{-1}(b, c, d)| = 6$ 
selects the set of polynomials with 3 distinct simple roots. These 
strata of $\C^3_{coef}$ are familiar under the names 
$\mathcal D_3$, $\mathcal D_2^\circ$ and $\mathcal D_1^\circ := \C^3 \setminus \mathcal D_2$. 

Note that, if a real cubic polynomial $P(z) = z^3 + bz^2 + cz + d$ has a 
single simple root, the 
triple $(b, c, d)$ is not in the image of the real Vi\`{e}te map.
\smallskip

The Jacobi matrix $D\mathcal V$ of the Vi\`{e}te map $\mathcal V$ is
\[\left( \begin{array}{ccc}
-1 & -1 & -1\\
v + w & w + u & u + v\\
-vw & -wu & -uv
\end{array}
\right)
\] 
and its determinant, the Jacobian $J\mathcal V$, is equal to $(v - u)(w - v)(u - w)$.
Therefore, away from the three planes $\Pi_{uv} := \{ v = u\}$, $\Pi_{vw} := \{ w = v\}$, 
$\Pi_{wu} := \{ u = w\}$  the 
rank of the Vi\`{e}te map is 3. On each the three planes it drops to 2, and at 
along the diagonal line $L := \{u = v = w\}$---to 1. Because of the $S_3$-symmetry, all the 
three planes have identical images under the $\mathcal V$. \smallskip

The Jacobian $J\mathcal V$ is not invariant under the permutations of the variables $u, v, w$. 
In fact, it changes sign under the transpositions of any two variables. Therefore, $J\mathcal V$ 
can \emph{not} be expressed in terms of the elementary symmetric polynomials in $u, v, w$, 
that is, in terms of the coefficients $b, c, d$. However, its square, \emph{the discriminant},  
\begin{eqnarray}
(J\mathcal V)^2 = [(v - u)(w - v)(u - w)]^2 
\end{eqnarray}
is invariant, and thus, is a polynomial $\Delta$ in $b, c, d$. A painful calculation 
(cf. [V]) shows that the discriminant
\begin{eqnarray}
\Delta(b, c, d) =  b^2c^2 - 4b^3d + 18bcd - 4c^3 - 27d^2. 
\end{eqnarray}
Under the $\mathcal V$, the equations $\{\Delta(b, c, d) = 0\}$ and 
$\{(v - u)(w - v)(u - w) = 0\}$ are equivalent. Evidently, the latter equation selects 
the case of multiple roots. In other words, the \emph{discriminant surface} of 
Section 3 can be defined by an equation of degree 4: 
\begin{eqnarray}
\mathcal D_2 :=   \{b^2c^2 - 4b^3d  + 18bcd - 4c^3 - 27d^2 = 0\}
\end{eqnarray}
Who could imagine from the first glance that this unpleasant formula 
hides such a nice geometry?  

Over $\C$, the surface $\mathcal D_2$ coincides with 
the $\mathcal V$-image of each of the planes $\Pi_{uv}$, $\Pi_{vw}$, $\Pi_{wu}$. 
Over $\R$, by a stroke of good luck, a similar conclusion holds: 
if a real cubic polynomial has a complex root of multiplicity $\geq 2$, then
all its roots must be real.

\begin{lem} The  surface $\mathcal D_2$ in (4.4) admits a $uv$-parameterization by
\begin{eqnarray}
(b, c, d) = ( -2u - v,\;  u^2 + 2uv,\,  -u^2 v).
\end{eqnarray}
\end{lem}
 
{\bf Proof.} Under the substitution (4.1),  
the equations (4.4) and $J\mathcal V = 0$
 are equivalent. Clearly, the second equation says 
that one of the roots must be of multiplicity $\geq 2$. Then, putting $u = w$, gives 
the desired parameterization of $\mathcal  D_2$ by $\mathcal V|_{\Pi_{vw}}$. 
Note that this restriction is an 1-to-1 map and onto, both over $\C$ and $\R$. 
Indeed, any permutation from $S_3$ or acts trivially on triples of the form $\{(u, v, u)\}$, 
or takes them to triples which do not belong to the plane $\Pi_{vw}$. 
\qed
\bigskip

The discriminant curve $\mathcal  D_3$  is the image 
of the diagonal line $L = \{u = v = w\}$ under the Vi\`{e}te map $\mathcal V$.

\begin{lem} The surface $\mathcal  D_2$ divides the space 
$\R^3_{coef}$ into two chambers 
$\mathcal U_{3}$ and $\mathcal U_{1}$, one of which 
represents real cubic polynomials with 3 real roots and the other --- 
with a single real root. The chamber $\mathcal U_{3}$ is characterized 
by the inequality $$\{b^2c^2 - 4b^3d - 4c^3 + 18bcd - 27d^2  > 0\}.$$
 
In turn, the curve $\mathcal  D_3$ divides $\mathcal D_2$ into two domains 
$\{\mathcal D_2^{\pm}\}$, one 
of which corresponds to the cubic polynomials of the form 
$(z - u)^2(z - v)$ with $u < v$, and the other --- with 
$u > v$.
\end{lem}

{\bf Proof.}\quad By its definition, the chamber $\mathcal U_{3}$ is the 
interior of the image  $\mathcal V(\R^3_{coef})$. The 
chamber $\mathcal U_{1}$ is the interior of the image of the set 
$\{(u, v, \overline{v}) \in \C^3_{root}\}$, where $u \in \R$, under the complex Vi\`{e}te map. 
Clearly, the two sets $\{(u, v, \overline{v}) \in \C^3_{root}\}$ and 
$\{(u, v, w) \in \R^3_{root}\}$ intersect along 
the set of real roots with one of the roots being of multiplicity 
$\geq 2$. By Lemma 4.1, the $\mathcal V$-image of those is the surface $\mathcal  D_2$. 

For  distinct real roots $u, v, w$, the discriminant $(J\mathcal V)^2 > 0$. 
At the same time, for a single real root $u$, \, 
$(J\mathcal V)^2 = [(v - u)(\overline{v} - v)(u - \overline{v})]^2 = 
[(v - u)(\overline{v} - u)]^2(\overline{v} - v)^2 < 0.$ 

The proof of the claim about the domains $\{\mathcal D_2^{\pm}\}$ is 
even simpler. \qed 
\bigskip

Many geometric properties of the discriminant curve and surface, established 
in Section 3, can be easily derived employing the Vi\`{e}te map. Here are a few 
examples. 

Let $u, v, w$ be the roots of 
a polynomial  $z^3 + bz^2 + cz + d$. Recall that $\mathcal  D_2 = \mathcal V(\Pi_{uw})$. 
This gives its $uv$-parameterization (4.5).
Putting $u = v$ in (4.5), generates a familiar parameterization 
$(b, c, d) = A(u) = (-3u, \, 3u^2, -u^3)$ of the discriminant 
curve $\mathcal D_3$. 
 
A $t$-parametric equation of a generic tangent line to the curve $\mathcal D_3$
can be written as $A(u, t) =  (-3u,\, 3u^2, -u^3) + t(-3,\, 6u, -3u^2)$. 
Hence, the formula 
$(b, c, d) =  (-3u - 3t,\,  3u^2 + 6ut,\,  -u^3 - 3u^2t)$
describes a ruled surface, which has to be compared with the discriminant 
surface $\mathcal  D_2$ parameterized by (4.5). In order to show that the 
two surfaces coincide, we need to solve for  $t$ the system of equations:
\begin{eqnarray}
-2u - v& = & -3u - 3t  \nonumber\\
u^2 + 2uv & = & 3u^2 + 6ut \nonumber\\
-u^2 v& = & -u^3 - 3u^2t.
\end{eqnarray}
The only solution is given by $t = (v - u)/3$, in other words, 
$A(u, \frac{1}{3}(v - u)) = \mathcal V (u, v)$. We have arrived to a 
familiar conclusion: $\mathcal D_2$ 
is a ruled surface comprised of lines, tangent to $\mathcal D_3$. Each 
of the lines is produced with the help of $A(u, t)$ by fixing a particular 
value of $u$ and varying $t$. Equivalently, it can be produced with 
the help of $\mathcal V(u, v, u)$ by fixing $u$ and varying $v$ (note that 
(4.6) are linear expressions in $t$ and $v$).
Since $\mathcal V: \Pi_{uw} \rightarrow \mathcal D_2$ is a 1-to-1 map, 
distinct  lines $\{u = u_\star\}$ in the $uv$-plane must have disjoint images
$T^{u_\star}_2 \subset \A^3_{coef}$.  
Therefore, for a given point $P \in \mathcal D_2$, there 
is a single line through $P$ and tangent to $\mathcal D_3$. 
\smallskip

Next, we will determine the equation of a generic plane $T$, tangent to 
the discriminant surface. At the point $\mathcal V(u, v, u)$, it is spanned by 
the two vectors $\partial_u \mathcal V(u, v, u) = (-2,\, 2u + 2v, -2uv)$ and  
$\partial_v \mathcal V(u, v, u) = (-1,\, 2u, -u^2)$. Unless $u =v$, the 
two tangent vectors are independent. As before, the vector $n(u) = (u^2, u, 1)$ is 
orthogonal to both vectors $\partial_u \mathcal V(u, v, u)$ and 
$\partial_v \mathcal V(u, v, u)$ and therefore, to $T = T(u,v)$. Furthermore,
since $n(u)$ is $v$-independent, the normal vector $n(u)$ is constant 
along the line  $T^u_2 = \{\mathcal V(u, v, u)\}_v \subset \mathcal D_2$. 
Since $T(u, v) \supset T^u_2$, it must be $v$-independent, and therefore, deserves 
the familiar name $T^u_1$. As before, the normal vector field $n(u)$ extends 
across the singularity $\mathcal V(u, u, u)$, and so is the distribution of tangent 
planes. 
\bigskip 

Let $w_\star$ be a fixed number. Consider the image of the plane 
$W^{w_\star} := \{ w = w_\star\}$ under the Vi\`{e}te map. Formulas (4.1) 
gives a $uv$-parametric description of that image: 
\begin{eqnarray}
(b, c, d) =  (-[u + v] - w_\star,\;\, [uv] + w_\star[u + v],\; -w_\star\cdot [uv]).
\end{eqnarray}
Solving  for $u + v$ and for $uv$, gives a 
very familiar linear relation 
$$w_\star^3 + bw_\star^2 + cw_\star + d = 0$$
among the variables $b, c, d$.
This leads to a still somewhat surprising conclusion: the image $\mathcal V(W^{w_\star})$ 
of the plane $W^{w_\star}$ is contained in the \emph {plane}\footnote{Note, that 
the $\mathcal V$-image of a \emph{generic} plane in $\A^3_{root}$ is a surface 
whose degree $> 1$.} $T^{w_\star}_1 \subset \A^3_{coef}$ tangent to the discriminant surface! 
\smallskip

The map $\mathcal V : W^{w_\star} \rightarrow T^{w_\star}_1$ is generically 2-to-1 map: 
the cyclic permutation group of order 2 acts on the triples of the 
form $\{(u, v, w_\star)\}_{u, v}$ by switching $u$ and $v$. The map is 1-to-1 along 
the diagonal line $\Delta^{w_\star} := \{(u, u, w_\star)\}$ in $W^{w_\star}$.

In the complex case, the map $\mathcal V : W^{w_\star} \rightarrow T^{w_\star}_1$ is 
\emph{onto} since it is algebraic and of the rank 2 at a generic point. 
In the real case, the semi-algebraic set 
$\mathcal V (W^{w_\star})$ occupies a region of the  plane $T^{w_\star}_1$, 
bounded by the curve $\mathcal V (\Delta^{w_\star})$. In fact, this curve, given 
by the parametric equation $(b, c, d) = (-2u - w_\star,\, u^2 + 2w_\star u,\, -w_\star u^2)$, is 
a parabola, which resonates with our experience with the \emph{quadratic} 
Vi\`{e}te map and its discriminant curve! Evidently, the region it bounds is the 
intersection of the chamber $\mathcal U_3 \subset \R^3$ (see Lemma 4.2)
with the plane $T^{w_\star}_1$. It can be characterized by the linear equation 
$\{w_\star^2b + w_\star c + d = -w_\star^3\}$ coupled with the quadric inequality 
$\{b^2c^2 - 4b^3d  + 18bcd - 4c^3  - 27d^2 > 0\}.$

All these observations are assembled in 

\begin{thm} The complex Vi\`{e}te map $\mathcal V$ takes each plane 
$W^{w_\star} := \{w = w_\star\}$ 
in $\C^3_{root}$ onto the plane $T^{w_\star}_1 := \{w_\star^2b + w_\star c + d = -w_\star^3\}$ 
in $\C^3_{coef}$, which is tangent to the discriminant surface $\mathcal D_2$. 
The map $\mathcal V : W^{w_\star} \rightarrow T^{w_\star}_1$  
is a 2-to-1 map, ramified along the quadratic curve $\mathcal V (\Delta^{w_\star}) 
\subset \mathcal  D_2 \cap T^{w_\star}_1$.
\smallskip

The real Vi\`{e}te map $\mathcal V$ takes each plane $W^{w_\star} \subset \R^3_{root}$  
onto the region of the tangent plane $T^{w_\star}_1$,
bounded by the parabola $\mathcal V (\Delta^{w_\star})$. As in the complex case, 
the map $\mathcal V : W^{w_\star} \rightarrow T^{w_\star}_1$  
is a 2-to-1 map, ramified along the parabola. \qed 
\end{thm} 

\begin{cor} The images of the three planes $W^{u}, W^{v} , W^{w}$
under the  Vi\`{e}te map are contained (in the complex case, coincide with) in the 
the planes tangent to $\mathcal D_2$ and passing through the point 
$P = \mathcal V(u, v, w)$. \qed 
\end{cor}

%%%%%%%%%%%%%%%%%%%%%%%%%%%%%%%%%%%%%%%%%%%%%%%%%%%

\section{A slice of reality and the reduction flow} 

In the search for formulas solving polynomial equations of degrees $d \leq 4$, 
the first step is to replace a generic equation by  
an equation of the \emph{reduced} form. The reduced form is based on polynomials 
of degree $d$ with no monomials of degree $d - 1$.  

The substitution $x = z - b/3$ transforms a generic cubic polynomial 
$P(z) = \break z^3 + bz^2 + cz + d$ to its reduced form 
$Q(x) = x^3 + px + q$. In this form the dimensions of the root and coefficient 
spaces are reduced by one. We can depict them using the comfortable 
geometry of the plane. 

Since, in the reduced case, the sum of the roots $u + v + w = 0$, we can 
take two of the roots, say $u$ and $v$, for the independent variables 
in the root plane. Then the reduced Vi\`{e}te map $\mathcal V$ can be written as
\begin{eqnarray} 
(c, d) =  (uv - [u + v]^2,\; uv[u + v]).
\end{eqnarray}

The discriminant in (4.3) collapses to
\begin{eqnarray}
[(v - u)(w - v)(u - w)]^2 =  -4d^3 -27c^2.
\end{eqnarray}

Motivated by the success of the substitution $z \rightarrow z - b/3$,  
we will examine the geometry of an 1-parametric group  of 
transformations $\{\Phi_t\}$ of the coefficient space, induced by the 
$t$-family of substitutions
$\{z \rightarrow z + t\}$. A typical transformation is described by the formula
\begin{eqnarray}
 \Phi_t(b, c, d) =  (b - 3t,\;\;  c - 2tb + 3t^2,\;\; d - tc + t^2b - t^3).
\end{eqnarray}

Formula (5.3) is the result of a straightforward computation of $P(z + t)$. 
By the Taylor formula, while 
$(b, c, d) = (\frac{1}{2}P''(0),\; P'(0),\; P(0))$, 
\begin{eqnarray}
\Phi_t(b, c, d) =  \big( \frac{1}{2}P''(-t),\; -P'(-t),\; P(-t)\big).
\end{eqnarray}

\begin{lem} The transformation $\Phi_t : \A^3_{coef} \rightarrow \A^3_{coef}$  
preserves the Hermitian or Euclidean volume in the coefficient space. 
\end{lem}

{\bf Proof.}\quad For each $t$, the linear part 
$(b, c, d) \rightarrow (b,\; c -2tb,\; d - tc + t^2b)$ of the affine transformation $\Phi_t$, 
defined by (5.3), has a lower-triangular matrix with the units along the diagonal. 
Thus, its determinant is equal to 1. \qed 
\smallskip

A  direct verification proves

\begin{lem}
For a fixed point $A = (b, c, d)$, the $t$-parametric curve $\Phi_t(b, c, d)$ is 
a solution of a  system of linear differential equations:
\begin{eqnarray}
\dot A(t) = \left( \begin{array}{rrr}
0 & 0 & 0\\
-2 & 0 & 0\\
0 & -1 & 0
\end{array} \right) 
A(t)\; + \; \left(
\begin{array}{r}
-3\\
0\\
0
\end{array} \right),
\end{eqnarray}
satisfying  the initial condition $A(0) = (b, c, d)^\ast$. \qed
\end{lem}

We call the flow $\{\Phi_t\}$ defined by (5.5) (equivalently, by (5.3) or (5.4)) the 
\emph{reduction flow}.

Denote by $\Psi_t: \A^3_{root} \rightarrow \A^3_{root}$ the translation  
by the vector $(t, t, t)$. 
Note that $\Psi_t$ and $\Phi_t$ are conjugated via the Vi\`{e}te map: 
$\mathcal V \circ \Psi_t = \Phi_t \circ \mathcal V$. Since 
the diagonal planes $\Pi_{uv}$, $\Pi_{vw}$, $\Pi_{wu}$ 
are obviously invariant 
under the $\Psi_t$-flow, their $\mathcal V$-images are invariant under the $\Phi_t$-flow. 
Therefore, employing Lemma 4.1, the $\Phi_t$-flow must preserve the discriminant surface 
$\mathcal D_2$ as well, as the discriminant curve $\mathcal D_3$. 

This observation can be reinforced. 
The function $(v - u)(w - v)(u - w)$ is clearly invariant under the flow $\Psi_t$. 
Therefore, the discriminant $\Delta(b, c, d)$ 
must be invariant under $\Phi_t$-flow. In particular, every surface of constant level 
of the polynomial $\Delta(b, c, d)$ is invariant 
under this flow. As a result, the whole web of  planes tangent to $\mathcal D_2$ and 
lines tangent to $\mathcal D_3$ is preserved under the transformation group $\{\Phi_t\}$: 
 each map $\Phi_t$ defined by (5.3) is an affine transformation. 

The proposition below captures these observations.

\begin{prop} The polynomial $\Delta(b, c, d) = b^2c^2 - 4b^3d + 18bcd - 4c^3 - 27d^2 $ 
is invariant under the 1-parametric group $\Phi_t$ defined by (5.3)---(5.5). 
In particular, the strata $\mathcal D_3 \subset  \mathcal D_2$ are $\Phi_t$-invariant. 
Therefore, the web of planes tangent to the discriminant surface, as well as the web 
of lines, tangent to the discriminant curve, is preserved by the $\Phi_t$-action.
\qed
\end{prop}

For a fixed number $k$, consider the plane $H^k \subset \A^3_{coef}$ defined by $\{b = k\}$. 
The reduced polynomials form the plane $H^0$.\smallskip 

We intend to slice the  ruled stratification
$\mathcal D_3 \subset \mathcal D_2 \subset \A^3$ by the planes 
$\{H^k\}$ and to investigate a typical slice together with its evolution under 
the flow $\{\Phi_t\}$.

We notice that each orbit $\{\Phi_t(P)\}_t$, where 
$P = (b, c, d)$, hits  the plane $H^k$ at a single point.  Indeed, (5.3) admits 
a single $t$ for which the $b$-coordinate of $\Phi_t(P)$ is $k$. Moreover, (5.3) 
implies that the orbit and the plane $H^k$ are transversal at the intersection. \smallskip

Consider a curve 
$\mathcal D_2^{[k]} = \mathcal D_2 \cap H^k$ and a point $\mathcal D_3^{[k]} = \mathcal D_3 \cap H^k$
in the slice $H^k$. The curve $\mathcal D_2^{[k]}$ is a cubic \emph{cusp}: just add the 
constraint $b = k$ to the parametrization (4.5) in Lemma 4.1. 

For any plane $T \subset \A^3_{coef}$, denote by $T^{[k]}$ its slice  $T \cap H^k$. 

Note that any plane $T^u_1$, tangent to $\mathcal D_2$, is in general position with the plane $H^k$: 
it contains the line $T^u_2$ (tangent to $\mathcal D_3$) which  is transversal to 
$H^k$. Thus, for any $u, k$, the intersection $T_1^{u [k]} = T^u_1 \cap H^k$ is a line. 
Furthermore,  this line $T_1^{u[k]}$ and the curve $\mathcal D_2^{[k]}$ must be \emph{tangent}: 
a transversal slice of two tangent surfaces produces a pair of tangent curves.  
Their point of tangency is the intersection of a line $T^u_2$, along which $T^u_1$ and 
$\mathcal D_2$ are tangent, with the slice $H^k$. 

In the complex case, through any point $P \in H^k \setminus \mathcal D_2^{[k]}$ 
there are exactly three tangent planes. Therefore, their $k$-slice consists of three lines 
which contain $P$ and are tangent to the discriminant cusp curve $\mathcal D_2^{[k]}$. Similarly, 
when $P \in \mathcal D_2^{[k]}$, there are two tangent planes through $P$, and their 
slice produces a pair of tangent lines to $\mathcal D_2^{[k]}$ (one of which is tangent at $P$). 
Finally, when $P \in \mathcal D_3^{[k]}$, the tangent plane through $P$ is unique. Its 
slice is the line in $H^k$ which contains the singularity $P$ of the cusp curve 
and is the limit of its tangents as they approach $P$.\smallskip 
 
The real case has a slightly different description. To visualize it, compare Figures 7 and 8. 

The real cusp $\mathcal D_2^{[k]}$ divides the plane $H^k$ into two regions 
$\mathcal U_3^{[ k]} := \mathcal U_3 \cap H^k$ and $\mathcal U_1^{[ k]} := \mathcal U_1 \cap H^k$.
There are three lines tangent to the curve $\mathcal D_2^{[k]}$ through every point in region 
$\mathcal U_3^{[ k]}$, and only one tangent line through every point of $\mathcal U_1^{[ k]}$.  
For any point on the cusp curve $\mathcal D_2^{\circ [k]}$, there are two tangent lines. Finally, 
through $\mathcal D_3^{[k]}$ there is a single line "tangent" to $\mathcal D_2^{[k]}$ at the apex. 

\begin{figure}[ht]
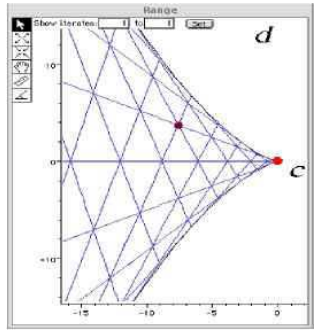

\centerline{\BoxedEPSF{cusp scaled 450}} 
\caption{The triangular, $S_3$-symmetric pattern is formed by the tangents to the discriminant curve
$\mathcal D_2^{[0]}$. Their {\it slopes differ by a fixed amount}. Each point in the domain $4c^3 + 27d^2 < 0$ is 
hit by tree tangent lines.}
\end{figure}

In short, the ruled stratified geometry of the slice is the slice of the ambient ruled 
stratified geometry. Moreover, the reduction flow respects these geometries.

\begin{prop} The section $\mathcal D_2^{[k]}$ of the discriminant surface $\mathcal D_2$ by 
the $cd$-plane $H^k := \{b = k\}$ is a $u$-parametric cubic curve $\{(c, d) = (-3u^2 - 2ku,\; 2u^3 + ku^2)\}$. 

The roots
of the polynomial $P(z) = z^3 + kz^2 + cz + d$ are equal to minus the slopes of the 
lines tangent to the curve $\mathcal D_2^{[k]}$ and passing through 
the point $P = (k, c, d)$.\smallskip 

The reduction flow $\Phi_t$ takes each curve $\mathcal D_2^{[k]}$ to the curve 
$\mathcal D_2^{[k - 3t]}$ and respects their webs of tangent lines. 
Specifically, if $u$ is a root of a polynomial $P(z) = \hfil\break z^3 + bz^2 + cz + d$,
then the  line $\{d = -uc - u^2b -u^3\}$ through $P = (b, c, d)$, residing in 
the plane $H^b$ and tangent to the curve $\mathcal D_2^{[b]}$, is mapped by $\Phi_t$ 
to the line \hfil\break $\{d = -(u + t)c - (u + t)^2b - (u + t)^3\}$ in $H^{b - 3t}$, 
passing through the point $\Phi_t(P)$ and tangent to the curve $\mathcal D_2^{[b - 3t]}$.
Thus, $\Phi_t$ is acting on the tangent lines by subtracting $t$ from their slopes. \qed
\end{prop}

\begin{cor} There is an invertible polynomial transformation $\mathcal K$ of the 
$bcd$-space $\A^3_{coef}$ mapping the discriminant surface $\mathcal D_2$ onto a surface 
$\tilde{\mathcal D}_2$ which is a Cartesian product of the cubic $\{4d^3 + 27c^2 = 0\}$ 
in the $cd$-plane and the $b$-axis $\A^1$. At the same time, $\mathcal K$ maps $\mathcal D_3$ 
onto the $b$-axis.
\end{cor}

{\bf Proof.}\quad Define a transformation $\mathcal K$ of the coefficient space 
as follows. This $\mathcal K$ moves any point $P = (b, c, d) \in H^b$ along its $\{\Phi_t\}$-trajectory 
until it arrives at a point $Q$ in the plane $H^0$; then it shifts $Q$ to a point $R$ 
in the plane $H^b$ with the same $cd$-coordinates as the ones of $Q$. The substitution 
of $t = b/3$ in the formula (5.3) helps to compute $\mathcal K(b, c, d)$ 
explicitly:
\begin{eqnarray}
\mathcal K(b, c, d) = (b', c', d') = (b,\;\;  c - \frac{5}{9}b^2, \;\; d - \frac{1}{3}bc - \frac{2}{27}b^3).
\end{eqnarray}
Because of its "upper triangular shape", the polynomial map $\mathcal K$ is invertible in 
the class of polynomial maps: one can uniquely express $(b, c, d)$ in terms of 
$(b', c', d')$. By Theorem 5.2, this  $\mathcal K$ has the desired properties. \qed

Thus, over the real numbers, there is a smooth homeomorphism of the pairs
$(\mathcal D_2 \subset \R^3) \approx (\R^2 \subset \R^3)$---the real surface 
$\mathcal D_2$ is topologically flat in the ambient space. 

%%%%%%%%%%%%%%%%%%%%%%%%%%

\section{stratified and ruled: the discriminants $\{\mathcal D_{d,k}\}$}

Notations in this section are similar, but more ornate than the corresponding
notations in Sections 2---5 (dealing with polynomials of degrees 2 and 3).\smallskip 

We consider the vector space $\A^d_{coef} = \mathcal D_{d,1}$ of monic polynomials 
$$P(z) = z^d + a_1z^{d-1} + ... + a_{d-1}z + a_d$$ 
of degree $d$ and 
its stratification  $\{\mathcal D_{d,k}\}_{1\leq k \leq d}$. Each strata $\mathcal D_{d,k}$
consists of polynomials with at least one of the roots being of multiplicity $\geq k$. 
Let $\mathcal D_{d,k}^\circ = \mathcal D_{d,k} \setminus \mathcal D_{d,k+1}$. 

As before, we can consider a subvariety $S_{d, k}$ in $\A^1  \times \mathcal D_{d,1}$,
defined by the system of equations
\begin{eqnarray}
\{P(z) = 0,\; P'(z) = 0,\; P''(z) = 0,\; ...\; ,\; P^{(k-1)}(z) = 0\}
\end{eqnarray}
Here $\A^1$ stands for the complex or real $z$-coordinate line and  
$P^{(j)}(z)$ denotes the $j$-th  derivative of $P(z)$.

Using the "upper triangular" pattern of (6.1), $\{a_d,a_{d-1}, ... , a_{d-k}\}$ can be 
uniquely expressed as polynomials  in $\{z, a_1, a_2, ... , a_{d-k-1}\}$. 
This produces a polynomial 1-to-1 parametrization $\mathcal H_{d,k}: \A^{d-k} \rightarrow 
S_{d, k}$ of the smooth variety $S_{d, k}$ of dimension $d-k$.

Evidently, that $\mathcal D_{d,k} \subset \mathcal D_{d,1}$ is the image 
of $S_{d, k}$ under the projection  
$\mathcal P: \A^1  \times \mathcal D_{d,1} \rightarrow \mathcal D_{d,1}$. 
As in the case of quadratic and cubic polynomials, for any $u \in \A^1$, 
the hyperplane $\{z = u\}$ hits $S_{d, k}$ along an $(d -k -1)$-dimensional \emph{affine space}
$N^u_{d,k}$  --- (6.1) are linear equations in the coefficients of $P(z)$. 
Hence, $S_{d, k}$ is a ruled variety. 

The projection $\mathcal P$ maps isomorphically $N^u_{d,k}$ onto an affine 
subspace $T^u_{d,k} \subset  \mathcal D_{d,k}$ of the coefficient space $\mathcal D_{d,1} $. 
This subspace parameterizes all monic  polynomials of degree $d$  having 
the root $u$ of multiplicity $\geq k$. Since any $P(z) \in \mathcal D_{d,k}$ 
has at least one root of multiplicity $\geq k$, the variety $\mathcal D_{d,k}$ 
is also comprised of the affine spaces $\{T^u_{d,k}\}$ of codimension one. However, they are not 
necessarily disjoint: a point in $\mathcal D_{d,k}$  can belong to many 
affine hypersurfaces. \smallskip

Each polynomial $P(z) \in \mathcal D_{d,k}$ has \emph{distinct} roots (real or complex) 
of multiplicities $\{\mu_1\geq\mu_2\geq ... \geq \mu_r\}$, with  $\mu_1 \geq k$ and $r \leq d$. 
Note that $P(z) \in \mathcal D_{d,k}^\circ$, iff $\mu_1 = k$.
In the complex case, the multiplicities $\{\mu_i\}$ define a \emph{partition} 
$\mu_P := \{\sum_{i=1}^r \mu_i = d\}$ 
of $d$. We also
interpret $\mu_P$ as a non-increasing function $\mu(i) = \mu_i$ on the set $\{1, 2, ... , d\}$, 
which takes non-negative integral values.
In the real case, the same interpretation holds, except 
that $\mu_P$ is a partition of a number which counts only the real roots with their
multiplicities. 

As a $\mu$-weighted configuration of roots 
deforms, a root of multiplicity $\mu_i$ and a root of multiplicity $\mu_j$ can 
merge producing a single root of multiplicity $\mu_i + \mu_j$. This results 
in a new partition $\mu'$ which we define to be \emph{smaller} than the original 
partition $\mu$. In a similar manner, several multiple roots can merge into a single one. 
Thus, the set of $d$-partitions acquires a \emph{partial ordering}: $\mu \succ \mu'$. 

For instance, the $d$-partition $\mu_{[k]}$, defined by the 
string of its values $(k, 1, 1, \; ...\; ,1)$, dominates any other $d$-partition
which starts with $k$.

In the real case,  complex roots occur in conjugate pairs of the same 
multiplicity, or are confined to the real number line $\R$. In what follows, 
while discussing the real case, we will use only the "$\R$-visible" part of 
the root configuration residing in $\R$. 
Only this part is captured by $\mu$. 
However, the partial ordering in the set of those $\mu$'s is induced in a way 
similar to the complex case. The only difference is that a pair of "invisible"
conjugate roots of a multiplicity $\mu_i$ can merge into a "visible" real root 
of multiplicity $2\mu_i$. 
For example, for $d = 4$, the "real"  
$\mu = (2, 0, 0, 0)$, corresponding to the configurations of one real root 
of multiplicity 2 and a pair of simple conjugate roots, is greater than 
the  real $\mu' = (2, 2, 0, 0)$, corresponding to the configurations 
of two real roots of multiplicity 2. 
\smallskip

By the definition of the projection $\mathcal P$, for any $P \in \mathcal D_{d,k}$, 
the cardinality of the preimage $\mathcal P^{-1}(P) \subset  S_{d, k}$ 
is the number of distinct roots of multiplicity $\geq k$ possessed by $P(z)$, 
that is, $|\mu_P^{-1}([k,d])|$. 
By the same token, the number of spaces 
$T^u_{d,k}$'s to which $P \in \mathcal D_{d,k}^\circ$ belongs is exactly 
$|\mu_P^{-1}(k)|$. At the same time, the number of hyperspaces $T^u_{d,1}$'s, passing
through $P \in \mathcal D_{d,1}$, is $|\mu_P|$---the cardinality of 
the support of the function $\mu_P$.
Since for $k > 1$, a \emph{generic} point $P \in \mathcal D_{d,k}^\circ$ 
corresponds to the partition $\mu_{[k]} = (k, 1, 1, \; ...\; ,1)$,\, 
$|\mu_P^{-1}(k)| = 1$ and there is a \emph{single} space $T^u_{d,k}$ passing 
through $P$. 
\smallskip 	

Because $\mathcal P: N^u_{d,k} \rightarrow T^u_{d,k}$ is an isomorphism, the differential $D\mathcal P$
of the projection $\mathcal P: S_{d, k} \rightarrow \mathcal D_{d,k}$ can only be of the ranks 
$d - k$  or $d - k -1$. The locus of points in $S_{d, k}$ where the rank drops is 
characterized by the property $D\mathcal P(\partial_z) = 0$. This happens when the 
gradients $\{\nabla_j\}$  of the $k$ functions in (6.1) defining $S_{d, k}$ are 
orthogonal to the vertical vector $\partial_z$. The $z$-component of the gradient 
vector $\{\nabla_j\}$ is exactly $P^{(j+1)}(z)$. Hence, the locus in question 
is characterized by a system as (6.1) with $k - 1$ being replaced by $k$. 
Therefore, it is the set $S_{d, k+1} \subset S_{d, k}$. Put $S_{d, k}^\circ := S_{d, k} \setminus S_{d, k+1}$.
As a result,  \emph{locally},  $\mathcal P: S_{d, k}^\circ \rightarrow 
\mathcal D_{d,k}^\circ $ is a smooth 1-to-1 map of maximal rank (an \emph{immersion}).
In particular, $\mathcal P: S_{d, 1}^\circ \rightarrow 
\mathcal D_{d,1}^\circ$ is a covering map with a fiber of the cardinality  $n$.
Hence, the singular locus of  $\mathcal D_{d,k}$ consists of $\mathcal D_{d,k+1}$ 
together with the self-intersections $\Sigma_{d,k}^\circ$ of  $\mathcal D_{d,k}^\circ$,
where each branch of $\mathcal D_{d,k}^\circ$ has a well-defined tangent space. 
In fact, $\Sigma_{d,k}^\circ$,  consists 
of polynomials $P(z) \in \mathcal D_{d,k}^\circ$, for which $|\mu_P^{-1}(k)| > 1$, 
\, $k > 1$. 

We notice that, for $k > d/2$, $|\mu_P^{-1}(k)| = 1$. Therefore, the 
immersion $\mathcal P: S_{d, k}^\circ \rightarrow \mathcal D_{d,k}^\circ$ 
is a regular \emph{embedding}, provided $k > d/2$. 
\begin{lem} For $k > d/2$,\; $\mathcal D_{d,k}^\circ$ is a smooth 
quasi-affine\footnote{that is, a Zariski-open set of an affine variety} 
subvariety of $\A^d_{coef}$. Hence, for $k > d/2$, the singular locus of 
$\mathcal D_{d,k}$ is $\mathcal D_{d,k+1}$. \qed
\end{lem}

\begin{exmp} $\mathbf{(d = 4).}$
\end{exmp} 
$\mathcal D_{4,1}^\circ \subset \C^4_{coef}$ is comprised of polynomials $P$
with the partition $\mu_P = (1, 1, 1, 1)$. The hypersurface $\mathcal D_{4,2}^\circ$
is comprised of polynomials with $\mu_P =  (2, 1, 1 ,0)$ or $(2, 2, 0, 0)$, and 
its self-crossing $\Sigma_{4,2}^\circ$---of polynomials with  $\mu_P =  (2, 2, 0, 0)$. 
The nonsingular surface $\mathcal D_{4,3}^\circ$
is comprised of polynomials with  $\mu_P =  (3, 1, 0 ,0)$. Finally, the smooth curve $\mathcal D_{4,4}$
corresponds to the partition $(4, 0, 0, 0)$. Thus, each point of $\mathcal D_{4,1}^\circ$ belongs 
to four hyperplanes of the type $T^u_{4,1}$; each  point of the space 
$\mathcal D_{4,2}^\circ \setminus \Sigma_{4,2}^\circ$ 
---to a single plane $T^u_{4,2}$ and each  point of the surface $\Sigma_{4,2}^\circ$ 
---to two planes $T^u_{4,2}$; each  point of the surface $\mathcal D_{4,3}^\circ$ belongs to 
a single line of the type $T^u_{4,3}$.\smallskip

The case of real degree 4 polynomials is more intricate. Figure 9 
shows a \emph{slice} of this stratification by a hypersurface $\{a_1 = 0\}$ of reduced 
quadric polynomials. 

The partitions $\mu_P$ which 
correspond to the three chambers of $\mathcal D_{4,1}^\circ \subset \R^4$ are: $(1, 1, 1, 1)$ 
(four distinct real roots---the "triangular" chamber in Figure 9),  $(1, 1, 0, 0)$ 
(two distinct simple real roots---the chamber "below" the surface in Figure 9) and 
$(0, 0, 0, 0)$ (no real roots---the chamber "above" the surface). The 
space $\mathcal D_{4,2}^\circ$ is comprised of four chambers. 
Three walls, bounding in $\mathcal D_{4,1}^\circ$ the chamber of four distinct real roots, 
all correspond to $\mu_P =  (2, 1, 1, 0)$. The three walls are distinguished by the three 
orderings in which a root of multiplicity 2 and two simple roots can be arranged on 
the number line $\R$.
The fourth chamber of the hypersurface $\mathcal D_{4,2}^\circ$ corresponds to 
$\mu_P =  (2, 0, 0, 0)$  (in Figure 9, the two wings which, behind the triangular tail, merge 
into a smooth surface). Points of the surface $\Sigma_{4,2}^\circ$ --- the transversal self-intersection 
 of $\mathcal D_{4,2}^\circ$ --- correspond to $\mu_P =  (2, 2, 0, 0)$. 
In Figure 9 they form the upper edge of the triangular chamber. 
Points of the surface $\mathcal D_{4,3}^\circ$ correspond to  $\mu_P =  (3, 1, 0, 0)$. 
In Figure 9 they form the two lower cuspidal edges of the triangular chamber.   
Finally, the curve $\mathcal D_{4,4}^\circ$ corresponds to $\mu_P =  (4, 0, 0, 0)$. In 
Figure 9 it is the apex of the tail. \qed

\begin{figure}[ht]
\centerline{\BoxedEPSF{D(4,2)reduced scaled 800}}
\caption{This swallow's tail is the section of the real $\mathcal D_{4,2}$ by the hypersurface $\{a_1 = 0\}$}
\end{figure} 

\begin{exmp} $\mathbf{(d = 5).}$ 
\end{exmp}
To give a taste of structures to come, Figure 10 depicts a stratification of the space $\C^5_{coef}$
by the 5-partitions $\{\mu_P\}$ (this time represented by the Young type tableau---the graphs of the 
5-partition functions). They form 
a partially ordered set with its elements decreasing from the left to the right. 
By definition, $\mu \succ \mu'$, if a tableau $\mu$ can 
be  built from a tableau $\mu'$ by moving a few blocks to the \emph{left}.
Each $\mu$ is indexing a quasi-affine variety $\mathcal D_\mu^\circ$ in $\C^5_{coef}$ 
formed by polynomials whose roots have the multiplicities prescribed by $\mu$. 
Its closure $\mathcal D_\mu$ is a union $\cup_{\mu'\preceq \mu}\; \mathcal
D_{\mu'}^\circ$. The $\mu$'s with hook shaped tableaus produce the familiar stratification $\{ \mathcal D_{5,k}\}$.
The dimension of $\mathcal D_\mu$ is the number of columns in the tableau $\mu$. 
We will revisit this example many times. \qed
\bigskip

\begin{figure}[ht]
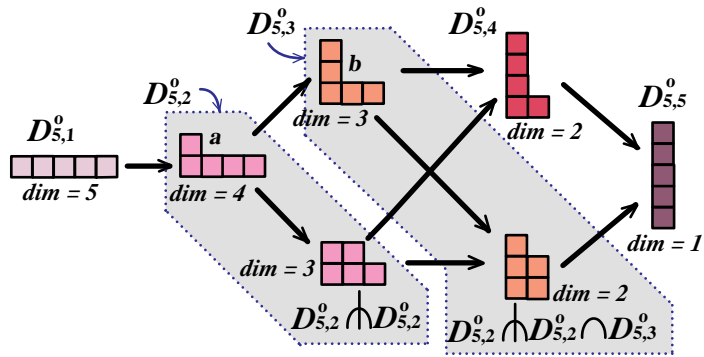

\centerline{\BoxedEPSF{deg5.tower scaled 500}}
\caption{The stratification of $\A^5_{coef}$ by the 5-partitions $\{\mu_P\}$.}
\end{figure}

We intend to show that each affine space $T^u_{d,k}$ is \emph{tangent} to the stratum 
$\mathcal D_{d,k+1}$ at some smooth point.
Recall that, the intersection multiplicity of this web of tangent spaces at a point  
$P \in \mathcal D_{d,k}^\circ$ is  $|\mu_P^{-1}(k)|$.
Moreover, we will see that the tangent cones to $\mathcal D_{d,k+1}$ span $\mathcal D_{d,k}$.
\smallskip

A normal space $\nu(S_{d,k})$ to the nonsingular variety 
$S_{d,k} \subset \A^1\times \mathcal D_{d,1}$ defined by (6.1), is spanned by $k$
independent gradient vectors $\{\nabla_j\}, 1 \leq j \leq k$. As (6.1) implies, 
at each point $(z, P) \in \A^1\times \mathcal D_{d,1} \approx \A^1\times\A^d$,
$$\nabla_j(z, P) = (P^{(j)}(z),\; n^{(j-1)}(z)),$$  
where $n^{(j-1)}(z)$ stands for the $(j - 1)$-st derivative of the vector 
\begin{eqnarray}
n(z) = (z^{d-1}, z^{d-2},\, ... ,\, z, 1).
\end{eqnarray}
However, at a point $(u, P) \in S_{d,k}$,\; $P^{(j)}(u) = 0$, provided $j < k$. 
Thus,
\begin{eqnarray} 
\nabla_j(z, P)& = &(0,\; n^{(j-1)}(z)),\; \mathrm{when}\; j < k; \nonumber\\
\nabla_k(z, P)& = &(P^{(k)}(z),\; n^{(k-1)}(z))
\end{eqnarray}

Therefore, the tangent space $\tau_x(S_{d,k})$ to $S_{d,k}$ at a point $x = (u, P)$ is 
spanned by vectors $w = (u,P) + (a, v) \in \A^1\times \A^d$, subject to constraints
\begin{eqnarray}
v \bullet n^{(j-1)}(u) &=& 0,\quad\quad 1 \leq j < k;\nonumber\\
v \bullet n^{(k-1)}(u) &=& - a \cdot P^{(k)}(u).
\end{eqnarray}
Here "$\bullet$" denotes the standard scalar product of vectors.  Evidently, if 
$v$ is orthogonal to all the vectors $\{n^{(j-1)}(u)\}_{1 \leq j < k}$, then it is possible 
to find an appropriate $a$ satisfying (6.6), provided $P^{(k)}(u) \neq 0$. When 
$P^{(k)}(u) = 0$, we have to consider the last equation from (6.6)  and 
to free the $a$.\smallskip

At the same time, the the $(d - k)$-dimensional space $T^u_{d, k-1}$ is defined by the 
linear equations
$$P(u) = 0,\; P^{(1)}(u) = 0,\; \; ...,\;\; P^{(k-2)}(u) = 0.$$ In terms of a vector 
$\tilde v \in \A^d$ they can be written as 
\begin{eqnarray}
\tilde v \bullet n^{(j-1)}(u) &=& -(u^d)^{(j-1)},\quad\quad 1 \leq j < k.
\end{eqnarray} 
Comparing the first $(k-1)$ equations from (6.4) with (6.5), we see that 
the first system of equations is the homogeneous part of the second system.
Therefore, the images of the space $N^u_{d, k-1}$ and the tangent space       
$\tau_{(u,P)}(S_{d,k})$ under the projection 
$\mathcal P: \A^1\times \mathcal D_{d,1} \rightarrow S_{d,1}$ coincide! Thus, 
$\mathcal P(\tau_{(u,P)}(S_{d,k})) = T^u_{d, k-1}$.\smallskip

The projection $\mathcal P: \A^1\times \mathcal D_{d,k} \rightarrow S_{d,k}$ takes 
the tangent space $\tau_x(S_{d,k})$ into the tangent cone $T_{k,P}$ of $\mathcal D_{d,k}$
at $P$. 
Since $\mathcal P: S_{d,k}^\circ \rightarrow \mathcal D_{d,k}^\circ$ is 
an immersion, the tangent cone $T_{k,P}$,\; $P \in \mathcal D_{d,k}^\circ$, is the 
union of the $\mathcal P$-images of the tangent spaces to $S_{d,k}^\circ$ at 
the points from $\mathcal P^{-1}(P)$. 
Therefore, for $P \in \mathcal D_{d, k}^\circ$, the cone $T_{k,P}$ is the union of 
$|\mu_P^{-1}(k)|$ affine spaces $\{T^u_{d, k-1}\}_{(u, P)}$, where $u$ runs 
over the set of distinct $P$-roots  of multiplicity $k$. 

The inclusion $ T^u_{d, k-1} \subset \mathcal D_{d, k-1}$ implies that, 
for $P \in \mathcal D_{d, k}^\circ$, the tangent cone $T_{k,P} \subset \mathcal D_{d, k-1}$. 
\smallskip

As an affine space, $\{T^u_{d, k-1}\}_{(u, P)}$ is determined by the equations 
$\{v \bullet n^{(j-1)}(u) \hfil\break = 0\}_{1 \leq j < k}$. Such a  set of equations
depends only on $u$, not on $P$. 
Therefore, it is shared by all the polynomials  
$P \in \mathcal D_{d, k}^\circ$ which have the \emph{same} root $u$ of multiplicity $k$.
In other words, along an open and dense set $T^u_{d,k}\cap \mathcal D_{d, k}^\circ$ in 
the $(d - k -1)$-space $T^u_{d,k}$, the tangent spaces 
$\{T^u_{d, k-1}\}_{P}$ 
are \emph{parallel} and, therefore, \emph{extend across (stabilize towards)
 the singularity} $T^u_{d,k} \cap \mathcal D_{d, k+1}$!
Furthermore, since $T^u_{d, k} \subset T^u_{d, k-1}$, as  affine spaces, all the  
 $\{T^u_{d, k-1}\}_{P \in T^u_{d, k}}$'s \emph{coincide}. 

For $k > d/2$, by Lemma 6.1, 
$\mathcal D_{d,k}^\circ$ is smooth, and the tangent bundle $\tau(\mathcal D_{d,k}^\circ)$
extends across the singularity $\mathcal D_{d,k+1} \subset \mathcal D_{d,k}$ to a vector bundle.
\bigskip

Although by now we understand the structure of the tangent cone $T_{k,P}$ at 
a generic point $P \in \mathcal D_{d, k}$ and the  stabilization of its 
components along some preferred directions towards the singular set 
$\mathcal D_{d, k+1}$, the structure of the tangent cone $T_{k,P}$ at 
at singular points $P \in \mathcal D_{d, k+1}$ still remains  uncertain. 
All what is clear that, for  $P \in \mathcal D_{d, k+1}^\circ$, the cone $T_{k,P}$
contains the well-understood tangent subcone $T_{k+1,P}$. \smallskip 

With this in mind, let's investigate in a more direct fashion the tangent cone  $T_{k,P}$ 
at a point $P(z) = (z - u)^k \hat P(z)$, where $\hat P(z)$ denotes a monic polynomial of 
degree $d - k$. When $P \in \mathcal D_{d,k}^\circ$,  $\hat P(u) \neq 0$. 

Let $P_t(z)$ be a smooth $t$-parametric curve in $\mathcal D_{d,k}$, emanating 
from the point $P(z)$. Locally, it can be written in the form $(z - u - a_t)^k [\hat P(z) + R_t(z)]$,
where $R_t(z)$ is a polynomial of degree $d-k-1$ and $lim_{t \rightarrow 0}\; a_t = 0$,  
$lim_{t\rightarrow 0}\; R_t(z) = 0$. The components of the velocity vector $\dot P_t$ to the $t$-parametrized
curve $P_t(z) \subset \mathcal D_{d,1}$ are the coefficients of the $z$-polynomial 
$$\dot P_t(z) = k(z - u - a_t)^{k-1}\, \dot a_t [\hat P(z) + R_t(z)] + (z - u - a_t)^k \dot R_t(z).$$ 
Since the curve $P_t(z)$ is  smooth at the origin,
$$\dot P_0(z) =  lim_{t \rightarrow 0} \dot P_t(z) = 
(z - u)^{k-1}[k \dot a_0 \hat P(z)  + (z - u) \dot R_0(z)]$$ 

Thus, a $\tau$-parametric equation of any \emph{line} from the tangent cone at $P(z)$ has a form
\begin{eqnarray}
P(z) + \tau \dot P_0(z) =\nonumber  \\ 
 (z - u)^k \hat P(z) + \tau 
(z - u)^{k-1}[ k \dot a_0 \hat P(z)  + (z - u) \dot R_0(z)]. 
\end{eqnarray}
As a $z$-polynomial, it is divisible by $(z - u)^{k-1}$---the tangent line resides in $\mathcal D_{d,k-1}$. 
Therefore, for \emph{any} $P \in \mathcal D_{d,k}$, $T_{k,P} \subset \mathcal D_{d,k-1}$.

Taking the (Zariski) closures, $\mathcal D_{d,k-1}$ contains the union 
of all tangent cones to $\mathcal D_{d,k}$.  On the other hand, since 
any $P \in \mathcal D_{d,k-1}$ is contained in some $T^u_{d,k-1}$ which 
is tangent to $\mathcal D_{d,k}^\circ$, we conclude that 
$\mathcal D_{d,k-1} = \tau(\mathcal D_{d,k})$ --- the union 
of all tangent cones to $\mathcal D_{d,k}$. 
\bigskip

For $k > d/2$, through each point $Q \in \mathcal D_{d,k}$ there is a single 
space $T^u_k$ tangent to $\mathcal D_{d,k+1}$. In particular, 
for any point $P \in \mathcal D_{d,k+1}$, there exist 
a single pair $T^u_{k+1} \subset T^u_k$ containing $P$. 
Consider an 1-dimensional space $L^u_k = L^P_k$ which contains 
$P$ and is orthogonal to $T^u_{k+1}$ in $T^u_k$.  We claim that the union
$\cup_{P \in \mathcal D_{d,k+1}}\, L^P_k = \mathcal D_{d,k}$. 
Furthermore, $\mathcal D_{d,k}$ is the space of a line bundle 
over $\mathcal D_{d,k+1}$ with a typical fiber $L^P_k$. Indeed,
any $Q \in \mathcal D_{d,k}$ belongs to a unique affine space 
$T^u_k \supset T^u_{k+1}$. Take the line in $T^u_k$ through $Q$
orthogonal to $T^u_{k+1}$. It hits $T^u_{k+1}$ at a point  
$P \in \mathcal D_{d,k+1}$. Thus, $Q \in L^P_k$. Since $k > d/2$, 
all the spaces $\{T^u_k\}$ are distinct and so are the lines 
$\{L^P_k\}$.
\smallskip

The preceding conclusions are summarized in the main result 
of this section---Theorem 6.1. In a way, it is a special case of our 
main result ---Theorem 7.1, but has a different flavor. Therefore, it 
is presented here for the benefit of the reader.

\begin{thm} Let $\A$ stand for the number field $\C$ or $\R$. 
Denote by $P \in \A^d_{coef}$ the point corresponding to 
a monic polynomial $P(z)$ of degree $d$.

\begin{itemize}
\item For any $1 \leq k \leq d$, the stratum $\mathcal D_{d,k}$ 
is a union of tangent cones to the stratum  $\mathcal D_{d,k+1}$.

\item Each stratum $\mathcal D_{d,k}^\circ \subset \A^d_{coef}$ is an 
immersed smooth manifold. For $k > d/2$, $\mathcal D_{d,k}^\circ$ is 
a smooth quasi-affine subvariety. Moreover, such a $\mathcal D_{d,k}$ 
is the space of a line bundle over $\mathcal D_{d,k+1}$.

\item Through each point $P \in \mathcal D_{d,k}^\circ$, there are exactly 
$|\mu_P^{-1}(k)|$ affine spaces $\{T^u_{d,k}\}_u$ tangent to the stratum 
$\mathcal D_{d,k+1}$. 
The spaces $T^u_{d,k}$ are indexed by the distinct $P(z)$-roots $\{u\}$ 
of multiplicity  $k$ over the  field $\A$. 
Each space $T^u_{d,k}$ is defined by the linear constraints 
$\{P(u) = 0,\; P^{(1)}(u) = 0,\; ... ,\; P^{(k-1)}(u) = 0\}$
imposed on the coefficients of $P(z)$. \hfil\break
For $k > 1$, through a generic point $P \in \mathcal D_{d,k}^\circ$ 
there is a single tangent space $T^u_{d,k}$. For $k > d/2$, every
point $P \in \mathcal D_{d,k}^\circ$ belongs to
 a single  space $T^u_{d,k}$.

\item Each space $T^u_{d,k}$ is tangent to $\mathcal D_{d,k+1}$ along the subspace $T^u_{d,k+1}$.

\item On the other side of the same coin, the tangent cone $T_{k,P}$ to $\mathcal D_{d,k}^\circ$ 
at a  point $P$ is the union of the affine spaces 
$\{T^u_{d,k-1}\}_u$, where $u$ is ranging over the distinct $P(z)$-roots of multiplicity $k$ 
over $\A$. \qed
\end{itemize}
\end{thm}
\smallskip

\begin{cor}  The problem of solving a polynomial equation $P(z) = 0$ over $\A$ is equivalent
to the problem of finding all hyperplanes $T$ passing through the corresponding point $P \in \A^d_{coef}$ 
and tangent\footnote{that is, belonging to a tangent space of a smooth point in $\mathcal D_{d,2}$.} 
to the discriminant variety $\mathcal D_{d,2}$. 

Specifically, 
consider the normal vector to such a hyperplane, normalized 
by the condition that its $d$-th component equals 1. 
Then, its $(d-1)$-st component gives a root $u$ of  $P(z)$. 
Via this construction, distinct roots $\{u\}$ of $P(z)$ over $\A$ and 
tangent hyperplanes $\{T\}$ through $P$ are in 1-to-1 correspondence. \qed 
\end{cor}
\bigskip

Let's return to Example 6.2 and Figure 10 to illustrate the claims of Theorem 6.1.
The open strata 
$\mathcal D_{5,1}^\circ,\, \mathcal D_{5,3}^\circ,\, \mathcal D_{5,4}^\circ,\, \mathcal D_{5,5}$ 
are smooth, while the stratum $\mathcal D_{5,2}^\circ$ has a transversal self-intersection
along a 3-fold $\mathcal D_{5,2}^\circ \cap \mathcal D_{5,2}^\circ$ which consists of points 
$P$ with $\mu_P = (2, 2, 1, 0, 0)$. The 3-dimensional strata $\mathcal D_{5,2} \cap \mathcal D_{5,2}$ 
and $\mathcal D_{5,3}$ are not in general position even at a generic intersection point: 
their intersection is a surface, not a curve. A generic point of 
$\mathcal D_{5,2} \cap \mathcal D_{5,2} \cap \mathcal D_{5,3}$ corresponds to the 
partition $\mu_P = (3, 2, 0, 0, 0)$. Similarly, the intersection of the surfaces $\mathcal D_{5,4}$ and 
$\mathcal D_{5,2} \cap \mathcal D_{5,2} \cap \mathcal D_{5,3}$ is the curve $\mathcal D_{5,5}$.
In short, the  more refined stratification $\{\mathcal D_\mu\}$ corresponding to the 5-partitions can be 
recovered from the geometry of the crude stratification 
$\mathcal D_{5,1} \supset \mathcal D_{5,2} \supset \mathcal D_{5,3}\supset \mathcal D_{5,4} 
\supset \mathcal D_{5,5}$.

There are 5 hyperplanes tangent to  $\mathcal D_{5,2}$ through every point of 
 $\mathcal D_{5,1}^\circ$, 4 hyperplanes through every point $P$ of 
 $\mathcal D_{5,2}^\circ$ with $\mu_P = (2, 1, 1, 1, 0)$ (that is, through every 
point of $\mathcal D_{5,2}^\circ \setminus (\mathcal D_{5,2}^\circ \cap \mathcal D_{5,2}^\circ)$),
3 hyperplanes through every point of  $\mathcal D_{5,3}^\circ$ with $\mu_P = (3, 1, 1, 0, 0)$ 
or through every point of $\mathcal D_{5,2}^\circ \cap \mathcal D_{5,2}^\circ$ with 
the $\mu_P = (2, 2, 1, 0, 0)$, 2 hyperplanes through every point of $\mathcal D_{5,4}^\circ$ 
or through every point in $\mathcal D_{5,3}$ with $\mu_P = (3, 2, 0, 0, 0)$, and finally, 1 
hyperplane through every point of $\mathcal D_{5,5}$. In short, the multiplicity of 
the web of tangent \emph{hyper}planes to the discriminant hypersurface at $P$ is the cardinality 
of the support of $\mu_P$ (which also happens to be the dimension of the 
stratum $\mathcal D_{\mu_P}$). 

At the same time, there is a single 3-space tangent to $\mathcal D_{5,3}$ through each 
point of $\mathcal D_{5,2}^\circ \setminus (\mathcal D_{5,2}^\circ \cap \mathcal D_{5,2}^\circ)$,
two 3-spaces  through each point of $\mathcal D_{5,2}^\circ \cap \mathcal D_{5,2}^\circ$,
a single plane tangent to $\mathcal D_{5,4}$ through each point of $\mathcal D_{5,3}^\circ$, 
and a single line tangent to $\mathcal D_{5,5}$ through each point of $\mathcal D_{5,4}^\circ$.
\qed  
\bigskip  
 
As in the discussion preceding Corollary 3.2, we can introduce the notions 
of a horizon and a projective horizon of a variety in $\A^d_{coef}$, as viewed 
from a point in its complement. Let's glance at a horizon of
the discriminant hypersurface in $\A^d_{coef}$.

\begin{cor} Over $\C$, any point $P \in \mathcal D_{d,1}^\circ$ has a 
horizon $hor(\mathcal D_{d,2}^\circ, P)$ comprised of $d$  codimension 2 affine 
spaces in $\C^d_{coef}$ which are in a general position. The variety $\mathcal D_{d,3}$ 
is inscribed in the horizon $hor(\mathcal D_{d,2}^\circ, P)$. Similarly, 
for any $P \in \mathcal D_{d,1}^\circ$, the projective horizon  
$hor_\pi(\mathcal D_{d,2}^\circ, P) \subset \P_P^{d-1}$ consists of 
$d$ hyperplanes in a general position. The variety $\pi_P(\mathcal D_{d,3})$ 
is inscribed in $hor_\pi(\mathcal D_{d,2}^\circ, P)$. \qed
\end{cor}

\bigskip

Now consider the smallest stratum --- $\mathcal D_{d,d}$. 
It is a smooth curve $\kappa (u)$ in $\mathcal D_{d,1}$ 
whose points correspond to polynomials of the form $(z - u)^d$. Its $u$-parametric 
representation is given by
\begin{eqnarray} 
\{ a_k = (-1)^k 
\left(
\begin{array}{l}
 d \\ k
\end{array} 
\right)
u^k \}_{1\leq k \leq d} 
\end{eqnarray}
With any point $\kappa (u)$ on the curve we associate a flag of vector 
subspaces $V^1_u \subset V^2_u  \subset \; ... \;  \subset V^{d-1}_u \;  \subset \A^d$ with 
the origins at $\kappa (u)$. 
Each osculating space $V^k_u$ is spanned by the linearly independent vectors 
$\kappa^{(1)}(u), \kappa^{(2)}(u), ... , \kappa^{(k)}(u)$ emanating from
$\kappa (u)$. Here 
$\kappa^{(j)}(u)$ stands for the $j$-th derivative of $\kappa (u)$ with 
respect to $u$. 

Remarkably, each vector $v \in V^{d-k}_u$, emanating from 
$\kappa (u)$, satisfies the first $k$ orthogonality conditions from (6.4), 
that is, 
\begin{eqnarray}
v \bullet n^{(j)}(u) &=& 0,\quad\quad 0 \leq j \leq  k-1.
\end{eqnarray}
In other words, as affine spaces, $V^{d-k}_u = T^u_{d,k}$\;!

In order to verify this claim, we have to check that $\kappa^{(q)}(u) \bullet n^{(p)}(u) = 0$
for each pair $(q, p)$, subject to $1 \leq q \leq d - k,\; 0 \leq p \leq k - 1$. 
The identity is a repackaging of the obvious 
identities $\partial^p_z\partial^q_u\{(z - u)^d\} |_{\{z = u\}} = 0$, being 
interpreted as scalar products of two vectors. Here $p + q \neq d$. 
\smallskip

These considerations, combined with Theorem 6.1, lead to Theorem 6.2 below. 
Similar statements can be found in [ACGH], pp. 136-137.
Theorem 6.2 testifies that all the geometric and combinatorial complexity 
of the discriminant varieties $\{\mathcal D_{d,k}\}_k$ can be derived from the geometry 
of a single curve $\mathcal D_{d,d} \subset \A^d$ !

\begin{thm} Each affine space, tangent to the variety $\mathcal D_{d,k}$, for an 
appropriate $u$, is of the form $V^{d-k+1}_u$. In different words, it is 
the $(d-k+1)$-th osculating space of the rational curve $\mathcal D_{d,d}$ at 
the point $\kappa(u)$.
Therefore, 
\begin{itemize}
\item the ruled variety $\mathcal D_{d,k}$ is the union of all 
osculating spaces $\{V^{d-k}_u\}_u$ at the points of the curve $\mathcal D_{d,d}$\, ;

\item the tangent cone $T_{k,P}$ to $\mathcal D_{d,k}$ at $P \in \mathcal D_{d,k}^\circ$
is the union of $|\mu^{-1}(k)|$ affine spaces $\{V^{d-k+1}_u\}_u$, 
where $u$ runs over the $P(z)$-roots of multiplicity $k$.\qed
\end{itemize}
\end{thm}

With Theorem 6.2 in place, the Fundamental Theorem of Algebra acquires a new 
geometric life. 

\begin{cor}{\bf (A geometrization of the Fundamental Theorem  \hfil\break of Algebra)}\hfil\break 
There exists an 1-to-1 algebraic map $\mathcal W: S^d(\mathcal D_{d,d}) \rightarrow 
\mathcal D_{d,1}$ from the $d$-th symmetric 
product of the complex rational curve $\mathcal D_{d,d} \subset \mathcal D_{d,1}$ onto 
the space $\mathcal D_{d,1}\approx \C^d$.  It is defined by the following geometric 
operation. 

For any point $X \in S^d(\mathcal D_{d,d})$, viewed as an unordered
collection of points $\{\kappa(u) \in \mathcal D_{d,1}\}$ with their  multiplicities 
$\{\mu_u\ \geq 1\}$\footnote{in other words, as an effective divisor $\sum_u \; \mu_u \kappa(u)$ of degree $d$  
on the curve $\mathcal D_{d,d}$},  
$\mathcal W(X)$ is defined to be the unique \emph{intersection point} 
of the osculating spaces $\{V^{d-\mu_u}_u\}_u$, taken at the points $\{\kappa(u)\}$.
\end{cor}

{\bf Proof.}\quad Any polynomial $P(z)$ is uniquely determined my its distinct roots $\{u\}$ 
with their multiplicities $\{\mu_u\}$. Therefore, there is a single point $P \in \mathcal D_{d,1}$
belonging to $\cap_u\; T_{d,\mu_u}^u$.  Moreover, any $P \in \mathcal D_{d,1}$ is such an 
intersection.  Now apply Theorem 6.2 which identifies $T_{d,\mu_u}^u$ with $V_u^{d-\mu_u}$.\qed
\smallskip
\smallskip

Let $P(V)$ denote the projective space associated with a finite dimensional 
vector space $V$. Denote by $V^\ast$ the dual vector space. Every hyperplane
$H \subset P(V)$ can be viewed a point $H^\vee \in P(V^\ast)$.  Recall, that the \emph{projective 
dual} $\mathcal D^\vee$ of a  variety $\mathcal D \subset P(V)$ is  
a subvariety of $P(V^\ast)$, defined as a closure of the set $\{H^\vee\}$ formed by the hyperplanes 
$H \subset P(V)$ \emph{tangent} to $\mathcal D$ at one of its smooth points. By definition, $H$ is 
tangent to $\mathcal D$ at a smooth point $x$, if $H$ contains the tangent 
space $\tau_x(\mathcal D)$ of $\mathcal D$ at $x$. For any projective variety $\mathcal D$, 
one has $(\mathcal D^\vee)^\vee = \mathcal D$. Furthermore, if $x$ is a smooth point of 
$\mathcal D$ and $H^\vee$---a smooth point of $\mathcal D^\vee$, then $H$ is tangent 
to $\mathcal D$ at $x$ if and only if $x$, regarded as a hyperplane $x^\vee$ in 
$P(V^\ast)$, is tangent to $\mathcal D^\vee$ at $H^\vee$ (cf. [GKZ], Theorem 1.1).

Generically, $\mathcal D^\vee$ is a hypersurface. When 
$codim(\mathcal D^\vee) > 1$, the original  $\mathcal D$ is a ruled variety. In the case of the 
special determinantal varieties $\{\mathcal D_{d,k}\}$, or rather their projectivizations, 
the dimensions of $\{\mathcal D_{d,k}^\vee\}$ drop drastically. Indeed, as Theorem 6.2 
implies, the tangent spaces of $\mathcal D_{d,k}$ at its smooth points form an 1-\emph{parametric}
family. For instance, $\mathcal D_{d,2}^\vee$ is just a curve! 
%%%%%%%
\begin{cor} 
\begin{itemize}
\item The dimension  $dim (\mathcal D_{d,k}^\vee) = k - 1$.  
\item Its degree $deg(\mathcal D_{d,k}^\vee) \leq  deg(\mathcal D_{d,k-1}) = (k - 1)(d - k + 2)$,
provided $k > 2$. 
\item  $\mathcal D_{d,2}^\vee$ is a curve 
$\mathcal C \subset P((\A^{d+1})^\ast)$ of degree $d$, and 
$\mathcal D_{d,2} = (\mathcal C)^\vee$. 
\end{itemize}
\end{cor}
{\bf Proof.}\quad We have seen that the family  $\{V^{d-k+1}_u\}_u$ of 
tangent affine spaces to $\mathcal D_{d,k} \subset \mathcal D_{d,1}$ is 1-dimensional. Each of these spaces
$V^{d-k+1}_u$  is contained in a $(k - 2)$-dimensional family $\mathcal H_u$
of hyperplanes. Hence, $k - 2 \leq dim (\mathcal D_{d,k}^\vee) \leq k - 1$. Since there are 
infinitely many spaces $\{V^{d-k+1}_v\}_v$ which  are not contained in any of the hyperplanes 
from $\mathcal H_u$, $dim (\mathcal D_{d,k}^\vee) = k - 1$.

The  $deg(\mathcal D_{d,k}^\vee)$ is the number of transversal intersection points of $\mathcal D_{d,k}^\vee$ 
with a generic affine space $W^{d - k + 1}$ contained in an affine chart of $\mathcal D_{d,1}^\vee$. 
Due to the projective duality, this number equals to the number of hyperplanes in $\mathcal D_{d,1}$ 
which contain a generic affine subspace $U^{k - 2} \subset \mathcal D_{d,1}$ and are tangent to 
$\mathcal D_{d,k}$. We can construct  $U^{k - 2}$ in a  way that links it with 
$\mathcal D_{d,k-1}$. 

Let $U^{k - 2} \subset \mathcal D_{d,1}$ be a generic affine subspace 
which hits $\mathcal D_{d,k-1}$ transversally at $deg(\mathcal D_{d,k-1})$ points $\{P_\alpha\}$. 
By a general position argument (the Bertini Theorem 8.18 in [H]), we can assume that all the $P_\alpha$'s
are smooth points in $\mathcal D_{d,k-1}^\circ$. For a smooth point $P_\alpha$,\; 
$\mu_{P_\alpha}^{-1}(k - 1) = 1$, provided $k > 2$. Therefore, in $\mathcal D_{d,k-1}$, there exists a single  
space $T^{u,\alpha}_{d,k-1}$ of dimension $d - k + 1$ tangent to $\mathcal D_{d,k}$ and passing 
through $P_\alpha$. Denote by $H_\alpha$ the minimal affine subspace in $\mathcal D_{d,1}$ which 
contains the transversal subspaces $T^{u,\alpha}_{d,k-1}$ and $U^{k-2}$ (whose intersection is $P_\alpha$). 
By its construction, $H_\alpha$ is a hyperplane which contains $U^{k-2}$ and is tangent to $\mathcal D_{d,k}$.  
In fact, any hyperplane $H$, which contains $U^{k-2}$ and is tangent to $\mathcal D_{d,k}$ at a point $P$,
 can be constructed in this way. Indeed, it must contain at least one of the spaces $T^{u}_{d,k-1}$ tangent 
to $\mathcal D_{d,k}$ at $P$. Because $U^{k-2}$ has been constructed in general position with 
$\mathcal D_{d,k-1} \supset T^{u}_{d,k-1}$, $U^{k-2}$ and $T^{u}_{d,k-1}$ must be in general position 
in $H$. Counting dimensions, $U^{k-2}$ and $T^{u}_{d,k-1}$ have a single point $P_\alpha$ of 
intersection. Therefore, $deg(\mathcal D_{d,k}^\vee) \leq deg(\mathcal D_{d,k-1})$. In order to replace 
the inequality by an equality, one needs to verify that all these tangent hyperspaces 
are distinct. We conjecture that this is the case.  By [Hi] (see also [W1], Theorem 2.2 and 
Corollary 2.3), $deg \mathcal D_{d,k} = k(d - k + 1)$. Therefore, $deg(\mathcal D_{d,k}^\vee) \leq 
(k - 1)(d - k + 2)$. In particular, $deg(\mathcal D_{d,3}^\vee) \leq  2d - 2$.

By a similar argument, the last claim of the theorem follows from the fact that a generic 
point of $\mathcal D_{d,1}$ is hit by $d$ hyperplanes tangent to $\mathcal D_{d,2}$. \qed 
\bigskip

As in Sections 2 and 4, the ruled geometry of the strata $\{ \mathcal D_{d,k} \}$ can be 
approached using the Vi\`{e}te map $\mathcal V_d : \A^d_{root} \rightarrow \A^d_{coef}$. 
It is defined by the elementary symmetric polynomials $\{\sigma_k(u_1, u_2, \; ...,\; u_d)\}$ 
which express the coefficient $a_k$ in terms of the roots $\{ u_i \}$.
The symmetric group $S_d$ acts on the space $\A^d_{root}$ by permuting the coordinates.
Denote by $St_U$ the stabilizer in $S_d$ of a point $U\in  \A^d_{root}$. 
If $U = (u_1, u_2, \; ...,\; u_d)$,  then 
as before, one can associate with $U$ a non-increasing function (a tableau) 
$\mu^U: \{1, 2, \; ..., \; r\} \rightarrow \{0, 1, \; ..., \; d\}$ which counts 
the numbers of equal coordinates in the string $U$. 
In these terms, $St_U \approx \prod_{i=1}^r S_{\mu^U_i}$, where $i$ runs over 
the support of $\mu^U$. 
Note that the cardinality of the preimage $\mathcal V_d^{-1}(\mathcal V_d(U))$ is the 
order $| S_d/ St_U | = d!/\prod_i  \{(\mu^U_i)!\}$. 

Distinct orbit-types $\{S_d/H\}_{H = St_U}$ give rise to 
a natural stratification $\{\A^{d, H^\circ}_{root}\}_H$ of the root space $\A^d_{root}$ and, 
because the Vi\`{e}te map is $S_d$-equivariant, --- to a familiar stratification 
$\{\mathcal D_{\mu_U}^\circ := \A^{d, H^\circ}_{coef}\}_H$ of the coefficient space 
$\A^d_{coef}$ (cf. Figure 10).
The coarse stratification $\{ \mathcal D_{d,k} \}$ can be assembled from this more 
refined stratification $\{\mathcal D_\mu^\circ\}_\mu$. In fact, over the complex numbers, 
$\mathcal D_{d,k}$ consists 
of all points $P = \mathcal V_d(U)$ for which, up to a conjugation, $St_U \supseteq S_k$.
Similarly, $\mathcal D_{d,k}^\circ$ is comprised  
of  $P = \mathcal V_d(U)$ for which, up to a conjugation, $St_U \supseteq S_k$ and 
does not contain any subgroup $S_j$ with $j > k$. Over the reals, the situation 
is more subtle: $ \mathcal V_d$ fails to be \emph{onto}. For instance, the  
Vi\`{e}te image of the hyperplane $\A^4_{root} \cap \{u_1 + u_2 + u_3 + u_4 = 0\}$ is 
not the whole space in Figure 9, but just the triangular chamber corresponding
to $\mu_P = (1, 1, 1 ,1)$.  

For a 
non-increasing function $\mu: \{1, 2, \; ..., \; r\} \rightarrow \{0, 1, \; ..., \; d\}$,  
so that $\sum_{q = 1}^r \mu_q \leq d$, denote by $K_\mu$ the vector subspace of $\A^d_{root}$ 
defined by the equations $\{u_i = u_j\}$, where  
$\sum_{q = 1}^p \mu_q \leq i, j \leq \sum_{q = 1}^{p+1} \mu_q$
and $p$ ranges over the support of $\mu$. The $\mathcal V_d$-image of 
this $K_\mu$  belongs to $\mathcal D_{d,\mu_1}$.  

For example, if $\mu = (4, 2, 2, 1)$, then  $K_\mu$ is defined by the equations 
$\{u_1 = u_2 = u_3 = u_4;\; u_5 = u_6;\; u_7 = u_8\}.$  The $\mathcal V_9$-image of 
this $K_\mu$ (of codimension 5) belongs to $\mathcal D_{9,4}$ and forms there a 
subvariety of codimension 2. In contrast, if $\mu$ corresponds to the 
partition $(4, 1,  1,  1,  1,  1)$, then  $\mathcal V_9(K_\mu)$ has codimension 0 
in $\mathcal D_{9,4}$.

Recall that $\mu_{[k]}$ is a partition, defined by $\mu_1 = k$ and, for $i > 1$, $\mu_i = 1$.
In the complex case, this $\mu_{[k]}$ describes a generic point of $\mathcal D_{d,k}$. 
Therefore, over the complex numbers,  $\mathcal V_d(K_{\mu_{[k]}}) = \mathcal D_{d,k}$.
Over the reals, simple complex roots generically occur in conjugate pairs, which allows 
for a greater variety of "generic"  $\mu$'s. In the previous example, in addition to 
the partition $9 = 4 + 1  + 1 + 1 + 1 + 1$, we must also consider "equally generic"  subpartitions 
$4 + 1  + 1 + 1 $ and $4 + 1$.
\smallskip

For any number $u$, denote by $\Pi^u_k$ the affine subspace of $\A^d_{root}$ defined by the equations 
$\{u_1 = u,\, u_2 = u,\, ... , \, u_k = u\}$. Evidently, $\mathcal V_d(\Pi^u_k) \subset \mathcal D_{d,k}$.
Furthermore, by the definitions, $\mathcal V_d(\Pi^u_k) \subseteq T^u_{d,k}$ and, over the complex numbers, 
$\mathcal V_d(\Pi^u_k) = T^u_{d,k}$. 

In view of the Theorem 6.1, we get the following proposition.

\begin{thm} Over $\C$, the Vi\`{e}te image of the vector subspace $K_{\mu_{[k]}}$ is 
the discriminant variety $\mathcal D_{d,k}$.
Any affine space $T$ of dimension $d - k + 1$, tangent to 
$\mathcal D_{d,k}$, is the image of some affine space $\Pi^u_{k-1} \subset \C^d_{root}$ under 
the Vi\`{e}te map $\mathcal V_d$.\smallskip  

Over $\R$, $\mathcal V_d(K_{\mu_{[k]}})$ forms a chamber in
the discriminant variety $\mathcal D_{d,k}$. Any affine space $T$ of dimension $d - k + 1$, 
tangent to $\mathcal D_{d,k}$, contains a chamber $\mathcal V_d(\Pi^u_{k -1})$. \qed
\end{thm}

As in Section 5, the product $[\prod_{i,j} (u_i - u_j)]$ is invariant under the 
the 1-parametric family of substitutions $\{u_k \rightarrow u_k + t\}$. As a 
result, the discriminant $\Delta_d(a_1, a_2,\; ... , \; a_d)$ is an invariant 
under the invertible algebraic transformations $\{\Phi_t: \A^d_{coef} \rightarrow \A^d_{coef}\}$ 
induced by the substitutions $\{z  \rightarrow z + t\}$. Hence, $\{\Phi_t\}$ preserve the 
hypersurface $\mathcal D_{d,2} \subset \A^d_{coef}$. Examining (6.1), we 
see that each variety $\mathcal D_{d,k}$ is invariant under the flow $\{\Phi_t\}$. 
In fact, each stratum  $\mathcal D_\mu^\circ$ is invariant as well --- the 
multiplicities of roots do not change under the substitutions $\{z  \rightarrow z + t\}$.
In particular, the curve $\mathcal D_{d,d}$ is a trajectory of the flow $\{\Phi_t\}$
which takes the point $(a_1, a_2,\; ... , \; a_d)$, representing a polynomial $P(z)$, 
to the point $(\tilde a_1, \tilde a_2,\; ... , \; \tilde a_d)$, where 
$\tilde a_k = \frac{(-1)^{d-k}}{(d-k)!} P^{(d-k)}(-t)$. For a fixed $t$, $\Phi_t$ is 
just an invertible linear transformation of $\A^d_{coef}$.  Expanding the formula for 
$\tilde a_k = \tilde a_k(t)$ as a polynomial in $t$, we see that $\{\tilde a_k(t) = a_k + b_k(t)\}$,
where $b_k(t)$ is a polynomial of degree $k$ with no free term. Because of the "upper triangular"
pattern of these formulas, the transformation $\Phi_t$ preserves the Euclidean volume form
of the space $\A^d_{coef}$.

Since each trajectory $\{\Phi_t(P)\}_t$ hits the hyperplane $\{a_1 = 0\}$ at a singleton $P_{red}$, 
the whole stratification $\{\mathcal D_\mu^\circ\}$ acquires a product structure:
$\mathcal D_\mu^\circ \approx \A^1 \times (\mathcal D_\mu^\circ)_{red}$, where 
$(\mathcal D_\mu^\circ)_{red} = \mathcal D_\mu^\circ \cap \{a_1 = 0\}$. In particular,
the swallow tail surface in Figure 9, being multiplied by $\R^1$, is isomorphic 
to the real discriminant 3-fold $\mathcal D_{4,2}$. 

These observations are captured in a well known lemma below which expresses the geometry 
of general polynomials in terms of the geometry of the reduced ones.
\begin{lem} The reduction flow $\{\Phi_t: \A^d_{coef} \rightarrow \A^d_{coef}\}$
preserves the stratification $\{\mathcal D_\mu^\circ\}$ of the coefficient 
space, as well as its Euclidean volume form. In particular, the webs of 
affine spaces, tangent to the strata $\mathcal D_{d,k}$, remain invariant under 
the flow. $\{\Phi_t\}$ also establishes the algebraic 
isomorphisms $\mathcal D_\mu^\circ \approx \A^1 \times (\mathcal D_\mu^\circ)_{red}$. 
\qed
\end{lem}
This completes our description of the stratification $\{\mathcal D_{d,k}\}$.

%%%%%%%%%%%%%%%%%%%%%%%%%%%%%%%%%%%%%%%%%

\section{the whole shebang: 
%stratified, but unruled
tangency and divisibility}

We are ready to extend results of the previous section to a generic stratum  $\mathcal D_\mu$.

By now, we have developed immunity to combinatorial 
complexities. This resistance will help us to meet the challenge of the $\mathcal D_\mu$'s 
intricate geometry. \smallskip 

Let $|aut(\mu)|$ denote the order the symmetry group of a partition $\mu = \{\mu_1 + \mu_2 + \; 
... \; + \mu_r\}$, i.e. 
all the permutations of the columns in the tableau $\mu$ which preserve its 
shape. Thus, $|aut(\mu)| = \prod_l[\#\mu^{-1}(l)]!$, where $l$ runs over the distinct
values of the function $\mu$. Put $|\mu| = r$.

\begin{lem} Each stratum $\mathcal D_\mu^\circ$ is a smooth quasiaffine variety in $\A^d_{coef}$ of 
dimension $|\mu|$. \footnote{Note  that,  $\mathcal D_{d, k}^\circ$, which can be singular, in general,  
consists of several $\mathcal D_\mu^\circ$'s.} Its degree 
$deg  (\mathcal D_\mu) = r ! \,\{\prod_{i = 1}^r \mu_i \}/ \{\prod_{l}[\#\mu^{-1}(l)]! \}$.
\end{lem}
{\bf Proof.} We generalize arguments centered on formulas (6.2)---(6.4). 

By definition, any polynomial from $\mathcal D_\mu^\circ$ is of the form 
$P(z) = \prod_{i = 1}^r\, (z - u_i)^{\mu_i}$ with all the roots $\{u_i\}$ being 
distinct. We can regard $\{u_i\}$ as coordinates in a space $\A^r$.  
Let $(\A^r)^\odot$ be an open subset of $\A^r$ --- the complement to 
the diagonal sets $\{u_i = u_j\}_{i \neq j}$.  

Denote by $S_\mu$ a subvariety 
of $\A^r \times \A^d_{coef}$ defined by the equations $\{P^{(j)}(u_i) = 0 \}$, where 
$1 \leq i \leq r$ and $0 \leq j < \mu_i$ (compare this with (6.1)). The equations 
claim that $u_1$ is a root of multiplicity $\mu_1$, $u_2$ is of multiplicity $\mu_2$, etc. 

Put $S_\mu^\odot := S_\mu \cap [(\A^r)^\odot \times \A^d]$.

The gradient of the function $P^{(j)}(u_i): \A^r \times \A^d_{coef} \rightarrow \A^1$ is given by 
the formula $\nabla_j(u_i, P) = (w_{i}^{(j)}(u_i),\;  n^{(j - 1)}(u_i))$.
Here the vector $w_{i}^{(j)}(u_i) \in \A^r$ has the number $P^{(j)}(u_i)$ as its 
$i$-th component, the rest of its coordinates vanish. The vector  $n^{(j - 1)}(u) \in \A^d$ 
is the $(j - 1)$-st derivative of the familiar vector $n(u) = (u^{d - 1}, u^{d - 2}, ... , u, 1)$.
At the points of $S_\mu^\odot$ the vectors $\{\nabla_j(u_i, P)\}$ are linearly 
independent; furthermore, their images $\{n^{(j - 1)}(u_i)\}$ under the projection 
$\mathcal P : \A^r \times \A^d_{coef} \rightarrow \A^d_{coef}$ are independent as 
well.
\smallskip 

By definition, the projection $\mathcal P$ 
takes $S_\mu^\odot := S_\mu \cap [(\A^r)^\odot \times \A^d]$ exactly onto 
$\mathcal D_\mu^\circ$.  It is an $|aut(\mu)|$-to-1 covering map: 
each polynomial in $\mathcal D_\mu^\circ$ determines the ordered list of 
its distinct roots $(u_1, ... , u_r)$ up to permutations from $aut(\mu)$.
Therefore, both $S_\mu^\odot$ and $\mathcal D_\mu^\circ$ are smooth quasiaffine
varieties. 
\bigskip

The degree $deg(\mathcal D_\mu)$ apparently has been computed by Hilbert in 
[Hi], however, I have to admit that I do not understand his arguments. A much more 
recent computation can be found in [W1], Theorem 2.2 and Corollary 2.3, pp. 377-78. 
There, using appropriate resolutions, the Hilbert function of $\mathcal D_\mu$ is 
calculated. 

Here is an alternative argument which I found easier to describe. Denote by 
$\vec{t}_k$ a $k$-vector  $(t, t, ... , t)$. Let $P_j^\mu(t_1, t_2, ... , t_r)$
denote the polynomial $\sigma_j(\vec{t}_{\mu_1}, \vec{t}_{\mu_2}, ... , \vec{t}_{\mu_r})$,
where $\sigma_j$ is the $j$-th elementary symmetric polynomial in $d$ variables.  
Evidently, $\{P_j^\mu(t_1, t_2, ... , t_r)\}_j$ define a parametrization $\mathcal V_\mu$ of 
$\mathcal D_\mu$. Therefore, the number of transversal intersections of a generic affine 
$(d - r)$-space with $\mathcal D_\mu^\circ$ is the number of solutions $(t_1, ... , t_r)$ of a generic 
linear system 
$\{\sum_{j=1}^r a_{ij}P_j^\mu(t_1, t_2, ... , t_r) = b_j\}_{1 \leq i \leq r}$, being divided 
by $\prod_l (\# \mu^{-1}(l))!$ --- the degree of the map $\mathcal V_\mu$ (and the order of 
the $\mu$-stabilizer).  Each of the polynomial 
$\sum_{j=1}^r a_{ij}P_j^\mu(t_1, t_2, ... , t_r) - b_j$ contains the same set of monomials
$\{t_1^{\nu_1}t_2^{\nu_2} ... t_r^{\nu_r}\}$, where $0 \leq \nu_i \leq \mu_i$. Therefore,
they all share the same Newton polytope of the volume $\mu_1\mu_2\ ...\mu_r$. By the 
Bernstein Theorem (cf. Theorem (5.4) in [CLO]), the number of solutions 
$(t_1, t_2, ... , t_r)$ ($t_i \neq 0$) of the generic system above is $r! [\mu_1\mu_2\ ...\mu_r]$. 
Without loss of generality, we can assume that all $t_i \neq 0$ --- otherwise, 
$P_r^\mu(t_1, t_2, ... , t_r)$ does not contribute its monomial to the Newton polytop.
Therefore, $deg(\mathcal D_\mu) = r! [\mu_1\mu_2\ ...\mu_r] / \prod_l (\# \mu^{-1}(l))!$.
\qed \bigskip 

Now, we need to introduce a few combinatorial notations. For any $d$-partition 
$\mu = (\mu_1, \mu_2,  ...\, , \mu_r)$
with $\mu_1 \geq \mu_2 \geq ... \geq \mu_r \geq 1$, denote by 
$\mu \downarrow \mathbf 1$ a $(d - r)$-partition defined by the sequence 
$(\mu_1 - 1,\, \mu_2 - 1,\,  ...\, ,\, \mu_r - 1)$.  It is 
supported on a smaller set of indices: $|\mu \downarrow \mathbf 1| = |\mu| - \#(\mu^{-1}(1))$. 
Also, let $\mu \uparrow \mathbf 1$ be a $(d - \#(\mu^{-1}(1)))$-partition defined by the 
rule: $(\mu \uparrow \mathbf 1)_i = \mu_i + 1$, when $\mu_i \geq 2$ and 
$(\mu \uparrow \mathbf 1)_i = 0$ otherwise.

Let $P(z) = \prod_i (z - u_i)^{\mu_i}$. Then  $P^{\downarrow \mathbf 1}(z)$ denotes 
the polynomial $\prod_i (z - u_i)^{\mu_i - 1}$ whose root multiplicities are described by 
$\mu \downarrow \mathbf 1$.
Also, let $\mathbf 1_r$ denote the partition $(1, 1, \,...\, , 1)$ of $r$.
\smallskip 

One can add partitions (cf. Figure 11): for a $d$-partition  $\mu = (\mu_1, \mu_2,  ...\, , \mu_r)$ 
and a  $d '$-partition  $\mu' = (\mu'_1, \mu'_2,  ...\, , \mu'_{r'})$, define 
a $(d + d')$-partition $\mu \uplus \mu'$ by the formula   
$$(\mu_1, \mu_2,  ...\, , \mu_r, \mu'_1, \mu'_2,  ...\, , \mu'_{r'}).$$  
Of course, the sequence above is no longer a monotone one. To get from it a Young-type 
tableau we need to reorder its terms.\smallskip

\begin{figure}[ht]
\centerline{\BoxedEPSF{sum scaled 180}}
\bigskip
\caption{}
\end{figure}

\begin{figure}[ht]
\centerline{\BoxedEPSF{gamma.count scaled 180}}
\bigskip
\caption{}
\end{figure}

Given a partition $\kappa$ of $m$ and a partition $\tau$ of $n$, $n \geq m$, 
we introduce $\gamma(\kappa, \tau)$ as the number of distinct monic polynomials 
of degree $m$, whose root multiplicities are dictated by $\kappa$, and which divide a 
particular monic polynomial of degree $n$ with the root multiplicities 
prescribed by $\tau$. In other words, $\gamma(\kappa, \tau)$ counts 
the number of different functions $\kappa': \{1, 2, ... , |\tau|\} \rightarrow \Z_+$,
such that: 1) for every $i$,  $\kappa_i' \leq \tau_i$ and 2) $\kappa'$ is the 
form $\sigma(\kappa)$, where $\sigma \in S_{|\tau|}$ is a permutation (cf. Figure 12).

A close formula computing $\gamma(\kappa, \tau)$ 
in terms of the $\kappa_i$'s and $\tau_j$'s is quite unappealing.  When 
$\gamma(\kappa, \tau) \neq 0$, we will write $\kappa \unlhd \tau$.
\bigskip

The proposition below summarizes most of what we know about the geometry 
of the varieties $\mathcal D_\mu^\circ$.
%%%%%%%%%%%%%%%
\begin{thm} {\bf ("Divisibility is tangency")}
\begin{itemize}
\item Over $\C$, for any point $P \in \mathcal D_\mu^\circ$ representing a 
polynomial $P(z) = \prod_i (z - u_i)^{\mu_i}$ 
\footnote{with all the $u_i$'s being distinct}, the tangent 
space $T_P\mathcal D_\mu^\circ \subset \C^d_{coef}$ to $\mathcal D_\mu^\circ$ 
at $P$ consists of all polynomials $Q(z)$ divisible by 
$P^{\downarrow \mathbf 1}(z) =  \prod_i (z - u_i)^{\mu_i - 1}$.
Thus, it is defined in $\C_{coef}^d$ by a system of linear constraints 
$\{Q^{(j)}(u_i) = 0\}$, where the $u_i$'s range over the \emph{multiple} roots 
of $P(z)$ and \; $0 \leq j \leq \mu_i - 2$. \smallskip

\item The space $T_P\mathcal D_\mu^\circ$
is \emph{tangent} to $\mathcal D_\mu^\circ$ along the $|\mu^{-1}(1)|$-dimensional 
affine space $V$ of polynomials divisible by $\prod_{\{i:\; \mu_i \geq 2\}} (z - u_i)^{\mu_i}$. 
In turn, $V$ is defined by linear equations $\{Q^{(j)}(u_i) = 0\}$, where the $u_i$'s 
range over the \emph{multiple} roots of $P(z)$ and $0 \leq j \leq \mu_i - 1$. \smallskip

If the gaps between distinct values of the function $\mu$ all are greater than 1,
then $T_P\mathcal D_\mu^\circ \cap \mathcal D_\mu = V$ --- the  space $T_P\mathcal D_\mu^\circ$
is tangent to the variety $\mathcal D_\mu^\circ$ at any point of their intersection.
\smallskip

\item The intersection 
$T_P\mathcal D_\mu^\circ \;\cap\; 
\mathcal D_{(\mu \downarrow \mathbf 1) \uplus (\mathbf 1_{|\mu|})}^\circ$ 
is open and dense in the affine space $T_P\mathcal D_\mu^\circ$. 
Hence,  $\overline{T}\mathcal D_\mu^\circ$ --- the closure of the union of all tangent
spaces $\{T_P\mathcal D_\mu^\circ\}_P$ --- coincides with   
$\mathcal D_{(\mu \downarrow \mathbf 1) \uplus (\mathbf 1_{|\mu|})}$. 
At the same time, 
$T\mathcal D_\mu^\circ = \coprod_{\{\mu': \; (\mu \downarrow \mathbf 1) \;\unlhd\; \mu' \}}
\mathcal D_{\mu'}^\circ$. \smallskip

\item Let a partition $\nu$ be such that $2\cdot \#(\nu^{-1}(1)) \geq |\nu|$. Then, 
reversing the flow in the previous bullet, 
$\mathcal D_\nu^\circ$ is an open and dense subset of $\overline T\mathcal D_\mu^\circ$,
where $\mu = \hfil\break (\nu \uparrow \mathbf 1) \uplus \mathbf 1_s$
with $s = 2\cdot \#(\nu^{-1}(1)) - |\nu|$.
Hence, such $\mathcal D_\nu$'s are ruled varieties. \smallskip

\item For any $d$-partitions $\mu, \mu'$ and each point $Q \in \mathcal D_{\mu'}^\circ$,
there are exactly \hfil\break $\gamma(\mu \downarrow \mathbf 1,\; \mu')$ $|\mu|$-dimensional 
spaces which are tangent to $\mathcal D_\mu^\circ$ and contain $Q$.
\end{itemize}
\end{thm}

{\bf Remark.} The theorem claims that a generic $\mathcal D_\mu^\circ$ has a very different geometry than the special
$\{\mathcal D_{d, k}^\circ \supset \mathcal D_{(k, 1, ... , 1)}^\circ\} $ indexed by the 
"hook-shaped" $\mu$'s.  For example, with
$d = 4$,  the tangent planes $T_P\mathcal D_{(2,2,0,0)}^\circ$  sweep 
$\mathcal D_{(1,1,1,1)} \setminus \mathcal D_{(4,0,0,0)}$ --- the surface $\mathcal D_{(2,2,0,0)}$ 
is "more curved" in $\C^4_{coef}$ than the surface $\mathcal D_{(3,1,0,0)}$ whose tangents span 
just a 3-fold (cf. Theorem 6.1). \smallskip 

{\bf Proof.} The argument is a refinement of arguments centered on formula (6.6).
For any $P \in \mathcal D_\mu^\circ$ representing a polynomial 
$P(z) = \prod_{i = 1}^r\, (z - u_i)^{\mu_i}$ with all its roots $\{u_i\}$ being 
distinct, consider a smooth $t$-parametrized curve in $\mathcal D_\mu^\circ$ emanating 
from $P$ and given by the formula: $P_t(z) = \prod_{i = 1}^r\, (z - u_i + a_i(t))^{\mu_i}$.
Let $\dot a_i := \frac{d}{dt} a_i(t)|_{t = 0}$ and  $\dot P(z) := \frac{\partial}{\partial t} P_t(z)|_{t = 0}$.
Then $\dot P(z) = \prod_{i = 1}^r\, (z - u_i)^{\mu_i} [\sum_{i = 1}^r\, \mu_i \dot a_i (z - u_i)^{-1}]$.
So, the $\tau$-parametric tangent line $P + \tau \dot P_t$ is represented by the polynomials  
\begin{eqnarray}
\prod_{i = 1}^r\, (z - u_i)^{\mu_i}[1 + \tau \sum_{i = 1}^r\, \mu_i \dot a_i (z - u_i)^{-1}]. 
\end{eqnarray}
They are divisible by $\prod_{i = 1}^r\, (z - u_i)^{\mu_i - 1}$ and are not divisible 
by $(z - u_i)^{\mu_i}$, unless $\dot a_i = 0$ or $\tau = 0$. Polynomials in (7.1) 
with all the $\dot a_i \neq 0$
correspond to  partitions of the form $\{(\mu \downarrow \mathbf 1) \uplus (\nu)\}_\nu$, with 
$(\mu \downarrow \mathbf 1) \uplus (\mathbf 1_{|\mu|})$ being the maximal element among them. 
When an $\dot a_i$ vanishes,  (7.1) becomes divisible by $(z - u_i)^{\mu_i}$. 

Let $\kappa: \{i\} \rightarrow \{0, 1\}$ be a book keeping function registering which $\dot a_i$'s 
vanish. Then (7.1) is divisible by $\prod_{i=1}^r (z - u_i)^{\mu_i - \kappa_i}$ and, 
\emph{generically}, is not divisible by 
$(z - u_i)^{\mu_i}$ for all $\kappa_i \neq 0$ or by $(z - u_i)^{\mu_i + 1}$ for all $\kappa_i = 0$. 
As a result, the corresponding tangent line $P + \tau \dot P_t$ is contained in 
$\mathcal D_{(\mu \downarrow \kappa) \uplus (\mathbf 1_{|\kappa|})} \subset 
\mathcal D_{(\mu \downarrow \mathbf 1) \uplus (\mathbf 1_{|\mu|})}$. Therefore, 
$T_P\mathcal D_\mu^\circ \subset \mathcal D_{(\mu \downarrow \mathbf 1) \uplus (\mathbf 1_{|\mu|})}$.
\smallskip

On the other hand, any monic polynomial $Q(z)$ of degree $d$, which is \emph{divisible} 
by $\prod_{i = 1}^r\, (z - u_i)^{\mu_i - 1}$, belongs to a tangent line as  in (7.1) 
from the tangent space of $\mathcal D_\mu^\circ$ at a point $P(z) = \prod_{i=1}^r (z - u_i)^{\mu_i}$.  
Here the \emph{simple} $P$-roots $\{u_i\}$ are chosen \emph{freely}. Indeed, put
$Q(z) = \prod_{i = 1}^r\, (z - u_i)^{\mu_i - 1}\hat Q(z)$, where $\hat Q(z)$ is monic of 
degree $|\mu|$. Comparing  $Q(z)$ with (7.1) leads to an equation 
\begin{eqnarray}
\hat Q(z) = 
\prod_{i = 1}^r\, (z - u_i)[1 + \tau \sum_{i = 1}^r\, \mu_i \dot a_i (z - u_i)^{-1}].
\end{eqnarray}  
This forces  $\{\dot a_i = \hat Q(u_i)/(\tau \mu_i \Delta_i)\}$, were 
$\Delta_i := \prod_{j \neq i} (u_i - u_j) \neq 0$.
With this choice of the velocity vector $\dot a$ at $P$, one gets an identity of 
monic polynomials of degree $|\mu|$,
which can be validated by comparing the LHS and RHS of (7.2) at $|\mu|$ distinct points $\{u_i\}$. 
Hence, $Q \in T_P\mathcal D_\mu^\circ$.

We notice that the tangent line $\{P + \tau \dot P\}$, $P(z) = \prod_i (z - u_i)^{\mu_i}$, can 
contain \emph{some} points $Q$ representing polynomials which are divisible by $(z - u_i)^{\mu_i + 1}$ 
or even by higher powers of $(z - u_i)$.  These $Q$'s are not in 
$\mathcal D_{(\mu \downarrow \kappa) \uplus (\mathbf 1_{|\kappa|})}^\circ$, but in its closure. 

We have shown that the affine space $T_P\mathcal D_\mu^\circ$ \emph{is comprised of polynomials 
divisible by} $P^{\downarrow \mathbf 1}(z)$. Therefore, $T_P\mathcal D_\mu^\circ$ is 
defined by linear equations $\{Q^{(j)}(u_i) = 0\}$, where the $u_i$'s range over the 
\emph{multiple} roots of $P(z)$ and \; $0 \leq j \leq \mu_i - 2$. Polynomials $Q(z)$
of degree $d$ which are divisible by $R(z) := \prod_{\{i:\; \mu_i \geq 2\}} (z - u_i)^{\mu_i}$ 
and have the rest of their roots simple, are clearly in 
$\mathcal D_\mu^\circ \cap T_P\mathcal D_\mu^\circ$. At the same time, using 
the previous description of spaces, tangent to $\mathcal D_\mu^\circ$, we get 
$T_P\mathcal D_\mu^\circ = T_Q\mathcal D_\mu^\circ$ (as affine spaces).
Therefore, $T_P\mathcal D_\mu^\circ$ is tangent to $\mathcal D_\mu^\circ$ along 
the subset of $\{Q(z)\}$ divisible by $R(z)$. Of course, not any polynomial
of degree $d$ which is divisible by $R(z)$ is in $\mathcal D_\mu^\circ$, 
but a generic one is. Note, that polynomials which are divisible by 
$R(z)$ form an affine space $V_R$ characterized by equations $\{Q^{(j)}(u_i) = 0\}$, 
where the $u_i$'s range over the \emph{multiple} roots of $P(z)$ and 
$0 \leq j \leq \mu_i - 1$. Hence, 
$V_R \cap \mathcal D_\mu^\circ \subseteq \mathcal D_\mu^\circ \cap T_P\mathcal D_\mu^\circ$ 
is the tangency locus of $\mathcal D_\mu^\circ$ and $T_P\mathcal D_\mu^\circ$.
It is open and dense in $V_R$. Both sets $\mathcal D_\mu^\circ \cap T_P\mathcal D_\mu^\circ$ and 
$\mathcal D_\mu^\circ \cap T_P\mathcal D_\mu^\circ$ are $|\mu^{-1}(1)|$-dimensional.

For example, when $P(z) = (z - 1)^3(z - 2)^2(z - 3)$, $P^{\downarrow \mathbf 1}(z) = 
(z - 1)^2(z - 2)$, and $T_P\mathcal D_\mu^\circ$ is comprised of polynomials of 
the form $(z - 1)^2(z - 2)(z - u)(z - v)(z - w)$. Here the roots $u, v, w$ are 
numbers of our choice. The polynomial $R(z) = \hfil\break (z - 1)^3(z - 2)^2$, and the line 
$V_R = \{(z - 1)^3(z - 2)^2(z - u)\}$ is contained in $\mathcal D_\mu$. 
Its intersection with $\mathcal D_\mu^\circ$ is characterized by the inequalities 
$u \neq 1, u \neq 2$. According to our argument, the 3-space $T_P\mathcal D_\mu^\circ$ 
is tangent to $\mathcal D_\mu^\circ$ along this line, pierced at two points. 
At the same time, $\mathcal D_\mu^\circ \cap T_P\mathcal D_\mu^\circ$ is a union 
of that line with the pierced line $\{(z - 1)^2(z - 2)^3(z - u)\}_{u \neq 1, 2}$ and 
the pierced curve $\{(z - 1)^2(z - 2)(z - u)^3\}_{u \neq 1, 2}$. Note that 
$T_P\mathcal D_\mu^\circ$ is not tangent to $\mathcal D_\mu^\circ$ along these two loci. 
\smallskip

For many partitions $\mu$, the complex intersection 
$\mathcal D_\mu^\circ \cap T_P\mathcal D_\mu^\circ$
simplifies to the linear form $\mathcal D_\mu^\circ \cap V_R$.
In particular, this happens when the gaps between distinct values 
of the function $\mu$ all are greater than 1 (i.e. the steps in the 
tableaux $\mu$ are higher than 1). For such a $\mu$, any 
polynomial from $\mathcal D_\mu^\circ$ which is divisible by 
$P^{\downarrow \mathbf 1}(z)$ is actually divisible by $R(z)$.
\smallskip

Recall, that each space tangent to $\mathcal D_\mu^\circ$  consists 
of polynomials which are divisible by $P^{\downarrow \mathbf 1}(z)$ 
for \emph{some} $P \in \mathcal D_\mu^\circ$.  
Therefore, in  order to prove the validity of the last statement of the 
theorem, we notice that the number $\gamma(\mu \downarrow \mathbf 1,\; \mu')$
measures exactly the number of distinct polynomials of the form 
$P^{\downarrow \mathbf 1}(z)$ which divide a given polynomial $Q(z)$,\, 
$Q \in \mathcal D_{\mu'}^\circ$. \qed\bigskip

We say that a $d$-partition $\mu$ is \emph{steep} if the gaps between the 
 distinct non-zero values of the (monotone) function $\mu: \{1, 2, ... ,  d\} 
\rightarrow \Z_+$ are greater than $|\mu|$. For example, a hook-shaped 
$\mu_{[k]}$ is steep when $2k > d + 2$.

\begin{cor} Assume that all $\mu_i \neq 2$. Then through each point $Q \in 
\mathcal D_{(\mu \downarrow \mathbf 1) \uplus (\mathbf 1_{|\mu|})}^\circ$ there 
is a unique $|\mu|$-space tangent to  $\mathcal D_\mu^\circ$.

For a steep $\mu$, all the tangent spaces to $\mathcal D_\mu^\circ$ are disjoint 
in $\C^d_{coef}$. Hence,  $T\mathcal D_\mu^\circ$ is a vector $[|\mu|- \#\mu^{-1}(1)]$-bundle 
over $\mathcal D_\mu^\circ$. 
\end{cor}

{\bf Proof.} Each $Q(z)$ as in the corollary has is a single divisor 
shaped by $\mu \downarrow \mathbf 1$. 
Similarly, for a steep $\mu$, 
any $d$-polynomial has no more than a single divisor shaped by 
$\mu \downarrow \mathbf 1$---the steps in the tableaux $\mu \downarrow \mathbf 1$ are 
"too tall". Hence, distinct tangent spaces are \emph{disjoint}. Each of them is 
tangent to $\mathcal D_\mu^\circ$ along a $\#\mu^{-1}(1)$-dimensional subspace. 
The orthogonal commplements to those (in the tangent spaces) provide the the bundle 
structure.
\qed 
\bigskip

Given a configuration of distinct points $\{x_j\}$ in the complex plane $\C^1$ (or 
in the complex projective space $\P^1$), equipped with  positive multiplicities 
$\nu_j$ \footnote{that is, an effective divisor $\sum_j \nu_j x_j$ of degree $d$.},
we define its \emph{resolution} to be a new collection points $\{y_{i,j}\}$ with 
 positive multiplicities $\mu_{i,j}$, so that $\nu_j = \sum_i \mu_{i, j}$ and 
$\{y_{i,j}\}_i$ reside in an $\epsilon$-neighborhood $U_j$ of $x_j$ \emph{free} 
of the rest 
of the points. As $\{\mu_{i, j}\}$ define a new partition $\mu$ of $d$,
the points $\{y_{i,j}\}_i$ "remember" their parent $x_j$. Specific locations 
of $\{y_{i,j}\}_i$ in $U_j$ are irrelevant, all we need is an association 
between $\{y_{i,j}\}_i$ and $x_j$ provided by $U_j$. This defines an equivalence 
relation between resolutions.

Now, fix $\mu \succ \nu$ and consider all equivalence classes of resolutions of 
$\sum_j \nu_j x_j$
for which $\{\mu_{i, j}\}_{i,j}$ produce $\mu$. We denote them 
$res(\sum_j \nu_j x_j,\; \mu)$, or alternatively, $res(P, \mu)$, where 
$P(z) = \prod_j (z - x_j)^{\nu_j}$. One can think of the set $res(P, \mu)$ as 
indexing \emph{locally} distinct branches of $\mathcal D_\mu^\circ$ in the 
vicinity of the point $P \in \mathcal D_\nu^\circ$. As long as 
$P \in \mathcal D_\nu^\circ$, all the sets $res(P, \mu)$ are isomorphic.

By eliminating all the $y_{i, j}$'s 
with $\mu_{i, j} = 1$ from our list, lowering the rest of $\mu_{i, j}$'s by 1, 
and still keeping the association of multiple $y_{i, j}$'s
 with $x_j$, analogous 
sets $res^{\downarrow\mathbf 1}(\sum_j \nu_j x_j,\; \mu) = 
res^{\downarrow\mathbf 1}(P, \; \mu)$ can be introduced. Again, the cardinality of 
$res^{\downarrow\mathbf 1}(P, \; \mu)$ depends only on $\mu$ and $\nu$.

\begin{cor} Let $\nu \prec \mu$ be two $d$-partitions. 
Let $\{P_t\}_{0 \leq t \leq 1}$ be a path in $\mathcal D_\mu$ 
so that, for $0 \leq t < 1,$ $P_t \in \mathcal D_\mu^\circ$ and 
$P_1 \in \mathcal D_\nu^\circ$. Then, as $t \rightarrow 1$, the tangent 
spaces $\{T_{P_t}\mathcal D_\mu^\circ\}$ \emph{stabilize} toward an affine 
$|\mu|$-dimensional space $T$ containing $P_1$. 

Although the limiting space $T$ can depend on the path
$P_t$ (which terminates at $P_1$), the number of such spaces at $P_1$ is finite. 
In fact, they are in an 1-to-1 correspondence with the elements of the set \,
 $res^{\downarrow\mathbf 1}(P_1, \; \mu)$.
In particular, when $\#(res^{\downarrow\mathbf 1}(P_1, \; \mu)) = 1$, the 
tangent bundle of $\mathcal D_\mu^\circ$ extends across the singularity 
$\mathcal D_\nu^\circ \subset \mathcal D_\mu$  to a vector bundle.
\end{cor}

For example, if $\nu = \{3 + 2 + 1\}$ and $\mu = \{2 + 2 + 1 + 1\}$, then the 
only $\mu$-resolution of $\nu$ is $\{(2 + 1) + 2 + 1\}$. Hence, 
$res^{\downarrow\mathbf 1}(P, \; \mu)$ consists of a single element 
$\mu\downarrow\mathbf 1 =\{ 2 + 2\}$, where the first 2 has 3 for a parent and
the second 2 has 2 for a parent. As a result, the tangent bundle to 
$\mathcal D_{(2,2,1,1,0,0)}^\circ$ extends across $\mathcal D_{(3,2,1,0,0,0)}^\circ$.

{\bf Proof.} Each space $T_{P_t}\mathcal D_\mu^\circ$ consists of polynomials 
divisible by $P_t^{\downarrow \mathbf 1}(z) =  \hfil\break \prod_i (z - u_i(t))^{\mu_i - 1}$.
As some of the distinct roots $\{u_i(t)\}$ merge when $t \rightarrow 1$, the set of 
polynomials divisible by $P_t^{\downarrow \mathbf 1}(z)$ converges to the set $T$ 
of polynomials divisible by a polynomial 
$Q(z) := lim_{t \rightarrow 1} P_t^{\downarrow \mathbf 1}(z)$ (also of degree 
$d - \#\{\mu^{-1}(1)\}$)\footnote{$Q(z)$ is different from $P_1^{\downarrow \mathbf 1}(z)$.}. 
Both sets are $|\mu|$-dimensional affine subspaces of $\C^d_{coef}$. 
Note that $P_1(z)$ is divisible by $Q(z)$.  Hence, $P_1 \in T$. 
Furthermore, the limiting polynomial $Q(z)$ does not depend on the choice of the path 
$P_t$, as long as the path is chosen so that the roots of $P_t(z)$, merging 
into a particular root of $P_1(z)$, are confined to its sufficiently small neighborhood
and their multiplicities are prescribed. 
In this context, "sufficiently small" means that the neighborhoods surrounding the 
roots of $P_1(z)$ are chosen to be disjoint. This prevents the roots of $P_t(z)$
from loosing focus on a parental root of $P_1(z)$ (that is, from "braiding" from a parental 
root to a different parental root). This remark justifies our previous definition of 
combinatorial resolution. Now it becomes clear that the limiting spaces at $P_1$ 
are in 1-to-1 correspondence with the elements of 
$res^{\downarrow\mathbf 1}(P_1, \; \mu)$, i.e. with the equivalence classes of multiple root
configurations governed by the $\mu$ and the association with  parental roots of $P_1$.  \qed
\smallskip
\begin{exmp}
\end{exmp}
Perhaps, an additional example can clarify Corollary 7.2 and the argument above. 
Take $\mu = \{3 + 2 + 1 + 1\}$ and $\nu = \{3 + 3 + 1\}$. Put $P_1(z) = 
(z - 4)^3(z - 6)^3(z - 8)$. We can resolve $P_1(z)$ only in two locally distinct ways: 
$R(z) = [(z - 3.9)^2(z - 4.1)]\times \hfil\break \times (z - 6)^3(z - 8)$ and  
$S(z) = (z - 4)^3[(z - 5.9)^2(z - 6.1)](z - 8)$,
each one  being consistent with the $\mu$. The tangent space to 
$\mathcal D_\mu^\circ$ at $R$ consists 
of monic polynomials of degree 7 which are divisible by 
$R^{\downarrow\mathbf 1}(z) = (z - 3.9)(z - 6)^2$ and the one at $S$ --- of 
polynomials divisible by $S^{\downarrow\mathbf 1}(z) = (z - 4)^2(z - 5.9)$.
The first is close to the limiting 4-space of polynomials divisible by 
$(z - 4)(z - 6)^2$, while the second is close to the limiting 4-space 
of polynomials divisible by $(z - 4)^2(z - 6)$. The two limiting 
spaces intersect along a 3-space of polynomials divisible by $(z - 4)^2(z - 6)^2$,
which happens to be the tangent space to $\mathcal D_\nu^\circ$ at $P_1$. \qed
\bigskip

The next proposition deals with the variety $\mathcal D_\mu^\vee$ projectively dual 
to $\mathcal D_\mu$. It resides  
in $\P(\C^d_{coef} \oplus \C^1)$. Unfortunately, for a general $\mu$, I do not know 
how to compute the degree of $\mathcal D_\mu^\vee$ (cf. Corollary 6.4).

\begin{cor} For any partition $\mu$,
$dim(\mathcal D_\mu^\vee) = d - 1 - \#\{\mu^{-1}(1)\}$
\end{cor}

{\bf Proof.} The arguments are similar to the ones in Corollary 6.4. 
By Theorem 7.1, $\mathcal D_{(\mu \downarrow \mathbf 1) \uplus (\mathbf 1_{|\mu|})}$
is $\bar{T}\mathcal D_\mu^\circ$. Therefore, $dim(\bar{T}\mathcal D_\mu^\circ) = 
dim(\mathcal D_{(\mu \downarrow \mathbf 1) \uplus(\mathbf 1_{|\mu|})}) = 
2|\mu| - \#\{\mu^{-1}(1)\}$. Since each tangent space is $|\mu|$-dimensional,
the whole family of these spaces must be $(|\mu| - \#\{\mu^{-1}(1)\})$-dimensional.
Each of the tangent spaces $T_P\mathcal D_\mu^\circ$ is contained in 
a $d - |\mu| - 1$ dimensional family of hyperplanes. Thus, 
$dim(\mathcal D_\mu^\vee) =  (|\mu| - \#\{\mu^{-1}(1)\}) + (d - |\mu| - 1) =
 d - 1 - \#\{\mu^{-1}(1)\}$.\qed 
\bigskip

The regular embedding $\mathcal D_\mu^\circ \subset \A^d_{coef}$ gives rise to a 
Gaussian map $G_\mu : \mathcal D_\mu^\circ \rightarrow PGr(|\mu|, d)$, where 
$PGr(k, d)$ denotes the Grassmanian of $k$-dimensional projective spaces in 
a $d$-dimensional projective space $\P(\A^d_{coef} \oplus \A^1)$. 
Let $\bar G_\mu(\mathcal D_\mu^\circ)$ stand for the closure of 
$G_\mu(\mathcal D_\mu^\circ)$ in $PGr(|\mu|, d)$. 

\begin{cor} The dimension of the variety $\bar G_\mu(\mathcal D_\mu^\circ)$ 
is $|\mu| - \#\{ \mu^{-1}(1)\}$, while the dimension of $\mathcal D_\mu^\circ$ 
is $|\mu|$. In particular, for $\mu = (k, 1, ... , 1)$,   
$dim(\bar G_\mu(\mathcal D_{d, k}^\circ)) = 1$, provided $k > 1$.
\end{cor}

{\bf Proof.} The union of tangent spaces to $\mathcal D_\mu^\circ$ spans an 
open set in the variety 
$\mathcal D_{(\mu \downarrow \mathbf 1) \uplus (\mathbf 1_{|\mu|})}$ of dimension 
$|(\mu \downarrow \mathbf 1) \uplus (\mathbf 1_{|\mu|})| = \#\{ \mu^{-1}([2, d])\} + |\mu|$.
Since the dimension of each tangent space is $|\mu|$, $dim(G_\mu(\mathcal D_\mu^\circ)) = 
\#\{ \mu^{-1}([2, d])\}$. \qed
\bigskip 

We can generalize the projective duality using Grassmanians instead of 
projective spaces. Given any projective 
variety $\mathcal X \subset \P(V^{n+1})$ of dimension $k$, denote by 
$\mathcal X^\vee$ the closure in 
$PGr(m, n)$, $k \leq m < n$, of all $m$-dimensional projective subspaces tangent 
(that is, containing a tangent space of $\mathcal X$) to $\mathcal X$ at its smooth points.

For an appropriate $s \in \Z_+$,  the number of $m$-dimensional 
projective spaces tangent to $\mathcal X$ and containing a fixed (but generic) $s$-dimensional 
projective subspace $U\subset \P(V^{n+1})$ is finite.
We define the \emph{degree} of $\mathcal X^\vee \subset PGr(m, n)$ to be this number. 

\begin{thm} For any partition $\mu$, consider the dual variety  
$\mathcal D_\mu^\vee$ residing in the Grassmanian 
$PGr(d -|\mu| + \#\{\mu^{-1}(1)\}, \; d)$. Then

\( \begin{array}{c}
deg(\mathcal D_\mu^\vee)\; \leq \;
\frac{(|\mu|\; + \;\#\{\mu^{-1}(2)\})!}{|\mu|!\; \cdot\; (\#\{\mu^{-1}(2)\})!}
\cdot deg(\mathcal D_{(\mu \downarrow \mathbf 1) \uplus (\mathbf 1_{|\mu|})}) \\ 
\;= \; \frac{(2|\mu|\; - \;\#\{\mu^{-1}(1)\}) !}{|\mu| !} \times \prod_{\{i:\; \mu_i > 2\}}
(\mu_i - 1) / \prod_{\{l \geq 2\}} (\#\{\mu^{-1}(l)\}!).
\end{array} \)\footnote{We conjecture that the estimate is sharp.}
\end{thm}

{\bf Proof.} From the definition, $deg(\mathcal D_\mu^\vee)$ is the number 
of affine subspaces 
of dimension $d -|\mu| + \#\{\mu^{-1}(1)\}$ in $\C^d_{coef}$ which are 
 tangent to $\mathcal D_\mu^\circ$ and containing a generic 
affine space $U$ of dimension $d - 2|\mu| + \#\{\mu^{-1}(1)\}$. 
Note that $dim(\mathcal D_{(\mu \downarrow \mathbf 1) \uplus (\mathbf 1_{|\mu|})})$
and $dim(U)$ are complementary. Pick $U$ to be transversal to
$\mathcal D_{(\mu \downarrow \mathbf 1) \uplus (\mathbf 1_{|\mu|})}^\circ$
at each of $deg(\mathcal D_{(\mu \downarrow \mathbf 1) \uplus (\mathbf 1_{|\mu|})})$
points $\{Q_\alpha\}$ of their intersection. Trough each point $Q_\alpha$, there are 
exactly  $\gamma(\mu \downarrow \mathbf 1, \; 
(\mu \downarrow \mathbf 1) \uplus (\mathbf 1_{|\mu|}))$ distinct spaces $T_{\alpha, j }$ 
which are tangent to $\mathcal D_\mu^\circ$. Therefore, there are at most 
$deg(\mathcal D_{(\mu \downarrow \mathbf 1) \uplus (\mathbf 1_{|\mu|})}) \times 
\gamma(\mu \downarrow \mathbf 1,\;  (\mu \downarrow \mathbf 1) \uplus (\mathbf 1_{|\mu|}))$ 
tangent $(d -|\mu| + \#\{\mu^{-1}(1)\})$-dimensional spaces which contain $U$. 
It remains to notice that 
$\gamma(\mu \downarrow \mathbf 1, \; (\mu \downarrow \mathbf 1) \uplus (\mathbf 1_{|\mu|}))$ 
is the number of choices of $\#\{\mu^{-1}(2)\}$ objects among 
$|\mu|\; + \;\#\{\mu^{-1}(2)\}$ objects. In particular, when 
$\#\{\mu^{-1}(2)\} = 0$, $deg(\mathcal D_\mu^\vee) \leq 
deg(\mathcal D_{(\mu \downarrow \mathbf 1) \uplus (\mathbf 1_{|\mu|})})$. 
This is the case for any $\mu = (k, 1, ..., 1)$ with $k > 2$. 

Finally, by [Hi], $deg(\mathcal D_{(\mu \downarrow \mathbf 1) \uplus (\mathbf 1_{|\mu|})}) = 
\theta \times \prod_{\{i:\; \mu_i > 2\}} (\mu_i - 1) / \prod_{\{l \geq 2\}} (\#\{\mu^{-1}(l)\} !)$,
where $\theta = \frac{(2|\mu|\; - \;\#\{\mu^{-1}(1)\})!}{(|\mu|\; + \;\#\{\mu^{-1}(2)\})!}$,
which completes the estimate for  $deg(\mathcal D_\mu^\vee)$. 
\qed
\bigskip

Although all the previous results were formulated for the affine or quasiaffine 
varieties $\mathcal D_\mu , \mathcal D_\mu^\circ$, in fact, many hold for 
their "projective  versions" $^\bullet\mathcal D_\mu , ^\bullet\mathcal D_\mu^\circ$. 
These are varieties of $\mu$-weighted configurations of (distinct) points in 
$\P^1$, in other words, the varieties of degree $d$ effective divisors $D$ of the form
$\sum_{i = 1}^{|\mu|} \mu_i p_i$, where $p_i = \hfil\break [a_i : b_i] \in \P^1$ are distinct points
 and $\{\mu_i > 0\}$.  
While the positively weighted configurations in $\A^1$ can be regarded as zeros of polynomials in 
one variable $z$, the positively weighted configurations in $\P^1$ can be regarded as zeros of 
homogeneous degree $d$ polynomials $P(z_0, z_1) = \prod_i \{-b_iz_0 + a_iz_1\}^{\mu_i}$ in two variables
$z_0, z_1$.

Adjusting the arguments  
 which have established Theorem 7.1 for a $t$-deformation 
$P_t(z_0, z_1) = \prod_i \{[-b_i - v_i(t)]z_0 + [a_i + u_i(t)]z_1\}^{\mu_i}$ 
of $P(z_0, z_1)$, we get

\begin{thm} 
\begin{itemize} \item
The space 
$T_D^\bullet\mathcal D_\mu^\circ$, tangent  to $^\bullet\mathcal D_\mu^\circ$ at a divisor $D$,  
consists of all the divisors $Q$ of the form $\sum_{i = 1}^{|\mu|} (\mu_i - 1) p_i + \hat Q$. Here 
$\hat Q$ is any effective divisor of degree $|\mu|$. 

\item In fact, $T_D^\bullet\mathcal D_\mu^\circ$ 
is tangent to $^\bullet\mathcal D_\mu^\circ$ along an open and dense set of a 
 $\#\{\mu^{-1}(1)\}$-dimensional projective space, 
formed by the divisors of the form $\sum_{\{i :\; \mu_i > 1\}} \mu_i p_i \hfil \break
 + \sum_{j = 1}^{\#\{\mu^{-1}(1)\}} q_j$,
where  $\{q_j\}$ are mutually distinct and distinct from the $p_i$'s.  

\item The tangent spaces $\{T_D^\bullet\mathcal D_\mu^\circ\}_D$ span a quasiprojective 
variety $T^\bullet\mathcal D_\mu^\circ$ formed by divisors of the form 
$\sum_{i = 1}^{|\mu|} (\mu_i - 1) p_i + \hat Q$, with $p_i$'s being distinct and $\hat Q > 0$ being 
of degree $|\mu|$.
\end{itemize}
\end{thm} 

Hence, for a given effective divisor $Q$ of degree $d$, the number of spaces 
$\{T_D^\bullet\mathcal D_\mu^\circ\}_D$ tangent to $^\bullet\mathcal D_\mu^\circ$ and containing 
$Q$ is the number of ways in which $Q = \sum_j \nu_j p_j$ can be represented as 
$\sum_{i = 1}^{|\mu|} (\mu_i - 1) p_i + \hat Q$, for some $\hat Q > 0$ and a collection of distinct $p_i$'s.
As we interpret partitions as non-increasing functions on the index set $\{1, 2, 3, ....\}$, they admit
another partial order: we say that $\nu \geq \mu$ when the function $\nu -\mu$ is non-negative 
(this partial order should not be confused with our old friend ---the partial order $\nu \succeq \mu$ induced 
by merging points in divisors).  Evidently, a divisor $Q = \sum \nu_j p_j$ belongs to the tangent space
$T^\bullet\mathcal D_\mu^\circ$, if and only if $\nu \geq \mu \downarrow \mathbf 1$. However, to compute the 
multiplicity of the tangent web $T^\bullet\mathcal D_\mu^\circ$ at $Q$ seems to be a tedious 
combinatorial problem: for given $\nu, \mu$, such that $\nu \geq \mu \downarrow \mathbf 1$, one needs 
to count the number of permutations $\sigma \in S_{|\nu|}$ which place the function 
$\sigma( \mu \downarrow \mathbf 1)$ \emph{below} $\nu$, divided by 
$\prod_l [\#\{(\mu \downarrow \mathbf 1)^{-1}(l)\} !]$---the order of the
stabilizer of $\mu \downarrow \mathbf 1$.  
\smallskip

{\bf Proof of Theorem 7.3.} First, we compute $\frac{d}{dt} P_t(z_0, z_1)|_{t = 0}$ which 
is given by the formula 
\begin{eqnarray}\prod_i (-b_iz_0 +a_iz_1)^{\mu_i - 1}\big [\sum_i \mu_i(-\dot v_iz_0 +\dot u_iz_1) 
\prod_{j \neq i} (-b_jz_0 + a_jz_1)\big ],
\end{eqnarray}
where $\dot u_i = \frac{d}{dt} u_i(t)|_{t = 0}, \dot v_i = \frac{d}{dt} v_i(t)|_{t = 0}$. This 
tells us that the tangent cone to $^\bullet\mathcal D_\mu^\circ$ at a divisor $D$ is contained 
in the set of divisors of the form $\sum_{i = 1}^{|\mu|} (\mu_i - 1) p_i + \hat Q$. Here  $\hat Q$ being 
an effective divisor of degree $|\mu|$. On the other hand, any homogeneous polynomial $Q(z_0, z_1)$ of degree
$d$ which is divisible by $\prod_i (-b_iz_0 +a_iz_1)^{\mu_i - 1}$ is of the form  (7.3) for an appropriate 
choice of the velocity vectors $\{\dot u_i, \dot v_i\}$. Indeed, let $Q(z_0, z_1) = \hat Q(z_0, z_1)
\prod_i (-b_iz_0 +a_iz_1)^{\mu_i - 1}$. Then from (7.3), the proportionality classes  $[\dot u_i : \dot v_i]$ 
are determined by the equations 
\begin{eqnarray} 
\hat Q(a_i, b_i) = \mu_i \left| \begin{array}{cc}
\dot u_i & \dot v_i\\
a_i & b_i
\end{array} \right| \cdot
\prod_{j \neq i}\left| \begin{array}{cc}
 a_j &  b_j\\
a_i & b_i
\end{array} \right|.
\end{eqnarray} 
In turn, for such a choice of the velocity vectors, $\sum_i \mu_i(-\dot v_iz_0 +\dot u_iz_1) \times \hfil\break 
\prod_{j \neq i} (-b_jz_0 + a_jz_1)$ is proportional to $\hat Q(z_0, z_1)$ (since, by (7.4), the two 
homogeneous polynomials of degree $|\mu|$ agree at $|\mu|$ distinct lines). \qed
\bigskip
%%%%%%%%%

Now we would like to make a few concluding remarks about  topology of the complex strata 
$\{\mathcal D_\mu\}$ and $\{\mathcal D_\mu^\circ\}$. A wonderfully rich account of 
the topological properties of disciminants, or rather their complements, can be found 
in [Va1], [Va2]. Both sources concentrate on more subtle description of real 
deteminantal varieties. A valuable topological information is also contained in [A1]---[A3],
[SW], [SK]. These papers tend to focus on calculations of the cohomologies of 
the complements to discriminant varieties of one kind or another.

First,  we notice that each $\mathcal D_\mu$ is a \emph{contractible} space. Indeed, 
the radial retraction of the plane $\C$ to its origin induces a retraction of any 
configuration to a singleton, taken with multiplicity $d$. 

%This homotopy reveals
%a (real) conic structure of $\mathcal D_\mu$. In fact, $\mathcal D_\mu$ 
%is a cone over a space which retracts onto the real variety 
%$\mathcal D_\mu^{S^1} \subset S^d(S^1)$. Its points are 
%$\nu$-weighted configurations in a unit circle $S^1$, for various 
%$\nu \preceq \mu$. The radial retraction of $\C^\ast$ on $S^1$ provides 
%a mechanism for the identification. 

\smallskip

In contrast, topology of  complex strata 
$\{\mathcal D_\mu^\circ\}$ and $\{\mathcal D_{d,k}^\circ\}$ is connected to the 
\emph{colored braid} groups similar to ones of Arnold [A1], [A2].  We think of distinct 
multiplicities of roots as being \emph{distinct colors}. 

Let $\mathcal U_k^\circ$ denote the configuration space of of $k$ \emph{ordered} distinct point 
in $\C$. Its fundamental group is the pure (or colored) braid group $F_k$ of $k$ strings.
By [FN], $\mathcal U_k^\circ$ is an Eilenberg-MacLane space $K(F_k, 1)$.

For a given partition $\mu : \{1, 2, ..., |\mu|\} \rightarrow \Z_+$, the space $\mathcal U_{|\mu|}^\circ$ 
is a finite covering of the space $\mathcal D_\mu^\circ$.  
The covering map identifies each ordered  configuration of $|\mu|$ distinct (simple) roots 
with a root configuration where roots acquire the multiplicities $\mu_i$'s and  roots of the same 
multiplicity do not enjoy  any order. Therefore, $\mathcal D_\mu^\circ$ must be an Eilenberg-MacLane 
space $K(\pi_1(\mathcal D_\mu^\circ), 1)$, where $\pi_1(\mathcal D_\mu^\circ)$ can be identified 
with  the $\mu$-colored braid group $B_{\mu}$. The number of strings in
the braids is $|\mu|$. A string starts  
and ends at roots of the same color-multiplicity $\mu_i$. In fact, $B_\mu$ is a subgroup 
of index $|\mu|!\;/ \prod_l [\#(\mu^{-1}(l))]!$ in $B_{|\mu|}$.

Due to Lemma 6.1, the slice of $\mathcal D_\mu^\circ$ by the hyperplane of reduced polynomials 
also is a $K(B_\mu, 1)$ space. \bigskip

Revisiting Example 6.2 and Figure 10, the space $\mathcal D_{5,1}^\circ$ 
is a $K(B_5, 1)$-space, where $B_5 := B_{(1, 1, 1, 1, 1)}$ is the braid group with 
5 strings. The stratum $\mathcal D_{(2, 1, 1, 1, 0)}^\circ$ is a 
$K(B_{(2, 1, 1, 1, 0)}, 1)$-space, where $B_{(2, 1, 1, 1, 0)}$ is the braid 
group of 4 strings with one of the strings colored with red and 3 strings with blue.
Similarly, both  $\mathcal D_{(3, 1, 1, 0, 0)}^\circ$ and $\mathcal D_{(2, 2, 1, 0, 0)}^\circ$ 
are $K(B_{(3, 1, 1, 0, 0)}, 1)$-spaces, where $B_{(3, 1, 1, 0, 0)}$ is the braid 
group of 3 strings with one of the strings colored with red and 2 strings with blue.
Both  $\mathcal D_{(4, 1, 0, 0, 0)}^\circ$ and $\mathcal D_{(3, 2, 0, 0, 0)}^\circ$ 
are $K(B_{(4, 1, 0, 0, 0)}, 1)$-spaces, where $B_{(4, 1, 0, 0, 0)}$ is a pure braid 
group of 2 strings. Finally, $\mathcal D_{(5, 0, 0, 0, 0)}$ is contractible. 
At the same time, $\mathcal D_{5,3}^\circ$, $\mathcal D_{5,4}^\circ$ both are 
$K(\Z, 1)$-spaces.\bigskip

I do not know whether, in general, $\mathcal D_{d,k}^\circ$ is an Eilenberg-MacLane space.
However, for $k > d/2$,   each 
$\mathcal D_{d,k}^\circ$ has a homotopy type of a circle. This is true because 
each $\mathcal D_{d,l}$, $l \geq k$, fibers over 
$\mathcal D_{d,l+1}$ with a fiber $\C$ and
each $\mathcal D_{d,l}^\circ$, $l \geq k$, fibers over 
$\mathcal D_{d,l+1}$ with a fiber $\C^\ast = \C \setminus \{0\}$. 
Thus, the fundamental group $\pi_1(\mathcal D_{d,k}^\circ) \approx \Z$, provided $k > d/2$. 
\bigskip

The space $\mathcal D_{d,1}^\circ$ admits an \emph{unordered} framing. Locally, 
it is comprised of $d$ independent vectors fields $\{n_k\}$. At a generic point 
$P(z) = \prod_k (z - u_k) \in \mathcal D_{d,1}^\circ$,\,  $n_k$ is the normal vector 
$(u_k^{d-1}, u_k^{d-2}, \; ... , \; u_k, 1)$ to 
the hyperplane $T_{d,1}^{u_k}$ which passes through $P$ and is tangent to the discriminant 
variety $\mathcal D_{d,2}$. Since all the roots are distinct, the Vandermonde matrix 
$Van(u_1, u_2,\; ...\; u_d)$, whose columns are the $\C$-linearly independent vectors
$\{n_k\}$, is of the maximal rank $d$. However, $\{n_k\}$ do not form vector fields:  
there is no consistent 
ordering attached to them. Nevertheless, they give rise to a well-defined algebraic embedding 
$$\mathcal Van^\circ: \mathcal D_{d,1}^\circ \rightarrow GL(d, \C)/S_d,$$
which induces a canonic epimorphism 
$$\mathcal Van^\circ_\ast : B(d,1) \rightarrow \pi_1(GL(d, \C)/S_d) \approx \Z \bowtie S_d.$$ 
Here the symmetric group $S_d$ acts on $\Z$ according to the parity of its permutations.

The embedding $\mathcal Van^\circ$ extends to an embedding  
$\mathcal Van: \mathcal D_{d,1} \rightarrow Mat(d, \C)/S_d,$ similarly defined in terms of the 
matrix $Van(u_1, u_2,\; ...\; u_d)$. 

The determinant gives rise to an obvious map 
$det: Mat(d, \C)/S_d \rightarrow \C/\{\pm 1\}$. 
In fact, $det[Van(u_1, u_2,\; ...\; u_d)] = \prod_{i < j}\, (u_i - u_j)$.
Evidently, $\mathcal D_{d,2} = 
(det \circ \mathcal Van)^{-1}(0)$, which also is the $\mathcal Van$-preimage of matrices of 
rank $d - 1$.  It is tempting to conjecture that the stratification $\{\mathcal D_{d,k}\}_k$ 
is a pull-back, under the embedding $\mathcal Van$, of the natural stratification 
$\{Mat_l(d, \C)\}_l$ 
of $(d \times d)$-matrices by rank $l$. However, the reality is different. Since the rank of 
a Vandermonde matrix $Van(u_1, u_2,\; ...\; u_d)$ is the number of distinct $u_i$'s,  
$\{\mathcal Van^{-1}(Mat_l(d, \C))\}_l$ defines a stratification 
$\{\mathcal R_{d,l}\}_l$ in $\C^d_{coef}$ by the number 
of \emph{distinct} roots, not by their maximal multiplicity as $\{\mathcal D_{d,k}\}_k$ does. 
In terms of the partitions $\mu$ associated with the roots,
$\mathcal D_{d,k}$ is comprised of polynomials $P$ with $max(\mu_P) \geq k$, while 
$\mathcal R_{d,l}$ is comprised of polynomials $P$ with $supp(\mu_P) \leq l$. 
For example, for $d = 4$, polynomials with $\mu = (3, 1, 0, 0)$ and with $\mu' = (2, 2, 0, 0)$
belong to $\mathcal R_{4,2}$. At the same time, polynomials with $\mu = (3, 1, 0, 0)$ belong 
to $\mathcal D_{4,3}$, while polynomials with $\mu' = (2, 2, 0, 0)$ belong to the larger 
stratum $\mathcal D_{4,2}$. In general, we only can claim that 
 $\mathcal D_{d,k} \subset \mathcal R_{d, d-k+1}$. \smallskip 
 
While the Vi\`{e}te map 
$\mathcal V: \C^d_{root} \rightarrow \C^d_{coef}$ transforms 
a simple linear stratification $\{\mathcal U_{d,k} := (u_1 = u_2 = ... = u_k)\}$ in $\C^d_{root}$ into
a "nonlinear" stratification $\{\mathcal D_{d,k}\}$ in $\C^d_{coef}$, the 
map $\mathcal Van$ has an "opposite" effect: it pulls back a linear stratification
$\{\mathcal Z_{d,k}\}$ in $Mat(d, \C)$ to produce $\{\mathcal D_{d,k}\}$.
A matrix $M = (m_{ij}) \in \mathcal Z_{d,k}$
when the first $k$ numbers among  $\{m_{d-1, j}\}_j$ are equal (if two such 
elements are equal, then so are the two columns of the Vandermonde matrix). 
In short, each stratum  $\mathcal U_{d,k}$ is a  linear subspace of $\C^d_{root}$, 
each stratum  $\mathcal Z_{d,k}$ is an linear subspace of  $Mat(d, \C)$.\smallskip

It is interesting to contemplate how do, with the help of $\mathcal Van^\circ$, the
topologies of $GL(d, \C)/S_d$ and  $\mathcal D_{d,1}^\circ$ interact. It seems that
 the induced cohomology  homomorphism 
$\mathcal Van^{\circ\ast}: H^\ast(GL(d, \C)/S_d)  \rightarrow H^\ast(B_d) = H^\ast(\mathcal D_{d,1}^\circ)$
is an epimorphism,  at least rationally. Perhaps, this interaction  is a subject for a different paper.

\end{document}